\documentclass[acmtoms]{acmtrans2m}

\usepackage{graphicx}
\usepackage{psfrag}
\usepackage{algorithm}
\usepackage{fancyvrb}
\usepackage{amsmath}
\usepackage{amssymb}
\usepackage[latin1]{inputenc}

\newcommand{\cgq}{$\mathrm{cG}(q)$}
\newcommand{\dgq}{$\mathrm{dG}(q)$}
\newcommand{\mcgq}{$\mathrm{mcG}(q)$}
\newcommand{\mdgq}{$\mathrm{mdG}(q)$}
\newcommand{\real}{\mathbb{R}}

\newcommand{\dt}{\, \mathrm{d}t}
\newcommand{\tab}{\hspace*{2em}}

\title{Algorithms and Data Structures for \\ Multi-Adaptive Time-Stepping}
\author{Johan Jansson \\ Royal Institute of Technology, Stockholm \and
        Anders Logg \\
        Center for Biomedical Computing, Simula Research Laboratory \\
        Department of Informatics, University of Oslo}
\date{\today}

\begin{document}

\begin{abstract}
  Multi-adaptive Galerkin methods are extensions of the standard
  continuous and discontinuous Galerkin methods for the numerical
  solution of initial value problems for ordinary or partial
  differential equations. In particular, the multi-adaptive methods
  allow individual and adaptive time steps to be used for different
  components or in different regions of space. We present algorithms
  for efficient multi-adaptive time-stepping, including the recursive
  construction of time slabs and adaptive time step selection. We
  also present data structures for efficient storage and interpolation
  of the multi-adaptive solution. The efficiency of the proposed
  algorithms and data structures is demonstrated for a series of
  benchmark problems.
\end{abstract}

\category{G.1.7}{Ordinary Differential Equations}{}[Error analysis, Initial value problems]
\category{G.1.8}{Partial Differential Equations}{}[Finite element methods]
\category{G.4}{Mathematical Software}{}[Algorithm design and analysis, Efficiency]

\terms{Algorithms, Performance}

\keywords{Multi-adaptivity, individual time steps, local time steps,
  multirate, ODE, continuous Galerkin, discontinuous Galerkin, mcgq, mdgq,
  C++, implementation, algorithms, DOLFIN}

\begin{bottomstuff}
  Johan Jansson,
  School of Computer Science and Communication,
  Royal Institute of Technology, SE--100 44 Stockholm, Sweden.
  \emph{Email}: \texttt{jjan@csc.kth.se}.
  \newline
  Anders Logg,
  Center for Biomedical Computing,
  Simula Research Laboratory,
  P.O.Box 134,
  1325 Lysaker, Norway.
  \emph{Email:} \texttt{logg@simula.no}.
  Logg is supported by an Outstanding Young Investigator
  grant from the Research Council of Norway, NFR 180450.
\end{bottomstuff}

\maketitle

\section{Introduction}

We have earlier in a sequence of papers
\cite{logg:article:01,logg:article:02,logg:article:08}
introduced the multi-adaptive Galerkin methods \mcgq{} and \mdgq{}
for the approximate (numerical) solution of ODEs of the form
\begin{equation}
    \begin{split}
      \dot{u}(t) &= f(u(t),t), \quad t\in(0,T], \\
      u(0) &= u_0,
    \end{split}
  \label{eq:u'=f}
\end{equation}
where $u : [0,T] \rightarrow \real^N$ is the solution to be computed,
$u_0 \in \real^N$ a given initial value,
$T>0$ a given final time,
and $f : \real^N \times (0,T] \rightarrow \real^N$ a given
function that is Lipschitz continuous in $u$ and bounded.

The multi-adaptive Galerkin methods~\mcgq{} and~\mdgq{} extend the
standard mono-adaptive continuous and discontinuous Galerkin
methods~\cgq{} and~\dgq{}, studied before in
\cite{Hul72a,Hul72b,Jam78,DelHag81,EriJoh85,Joh88,EriJoh91,EriJoh95a,EriJoh95b,EriJoh95c,EriJoh98,EriEst95,Est95,EstFre94,EstLar00,EstWil96,EstStu02},
by allowing individual time step sequences $k_i = k_i(t)$ for the
different components $U_i = U_i(t)$, $i = 1,2,\ldots,N$, of the
approximate solution~$U \approx u$ of the initial value
problem~(\ref{eq:u'=f}). For related work on local time-stepping, see
also~\cite{HugLev83a,HugLev83b,Mak92,DavDub97,AleAgn98,OshSan83,FlaLoy97,DawKir01,LewMar03,MR1485814,SavHun05,Sav08}. In
comparison with existing method for local time-stepping, the main
advantage of the multi-adaptive Galerkin methods~\mcgq{} and~\mdgq{}
is the automatic local step size selection based on a global a
posteriori error estimate built into these methods.

In the current paper, we discuss important aspects of the
implementation of multi-adaptive Galerkin methods. While earlier
results on multi-adaptive time-stepping presented
in~\cite{logg:article:01,logg:article:02,logg:article:08} include the
formulation of the methods, a priori and a posteriori error
estimates, together with a proof-of-concept implementation and results
for a number of model problems, the current paper addresses the
important issue of efficiently implementing the multi-adaptive methods
with minimal overhead as compared to standard mono-adaptive
solvers. For many problems, in particular when the propagation of the
solution is local in space and time, the potential speedup of
multi-adaptivity is large, but the actual speedup may be far from the
ideal speedup if the overhead of the more complex implementation is
significant.

\subsection{Implementation}

The algorithms presented in this paper are implemented by the
multi-adaptive ODE-solver available in
DOLFIN~\cite{www:dolfin,logg:preprint:06}, Dynamic Object-oriented
Library for FINite element computation, which is the C++ interface of
the new open-source software project
FEniCS~\cite{www:fenics,logg:article:12,logg:preprint:10} for the automation of
Computational Mathematical Modeling (CMM). The multi-adaptive solver
in DOLFIN is based on the original implementation Tanganyika,
presented in \cite{logg:article:02}, but has been completely rewritten
for DOLFIN and is actively developed by the authors.

\subsection{Obtaining the software}

DOLFIN is licensed under the GNU (Lesser) General Public License
\cite{www:LGPL}, which means that anyone is free to use or modify the
software, provided these rights are preserved.  The complete source
code of DOLFIN, including numerous example programs, is available at
the DOLFIN web page~\cite{www:dolfin}.

\subsection{Notation}

\begin{figure}[htbp]
  \begin{center}
    \psfrag{0}{$0$}
    \psfrag{i}{$i$}
    \psfrag{k}{$k_{ij}$}
    \psfrag{K}{$K_n$}
    \psfrag{T}{$T$}
    \psfrag{I}{$I_{ij}$}
    \psfrag{t1}{$t_{i,j-1}$}
    \psfrag{t2}{$t_{ij}$}
    \psfrag{T1}{$T_{n-1}$}
    \psfrag{T2}{$T_n$}
    \psfrag{t}{$t$}
    \includegraphics[width=12cm]{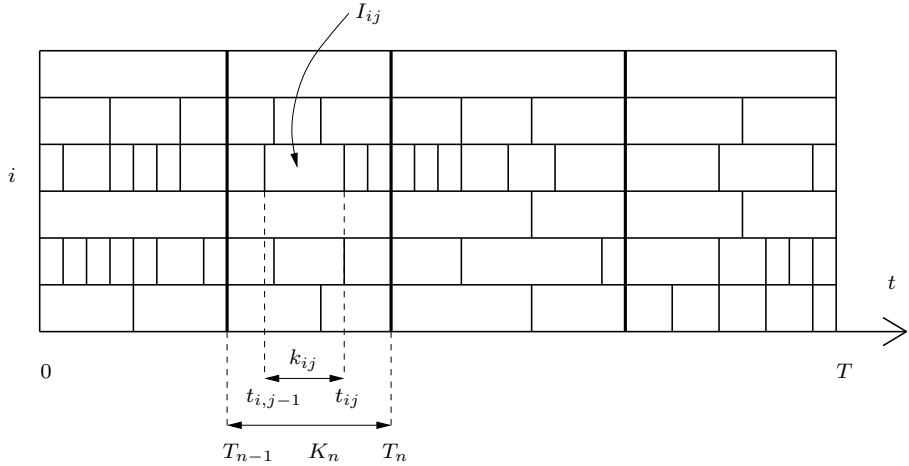}
    \caption{Individual partitions of the interval $(0,T]$ for different components. Elements
    between common synchronized time levels are organized
    in time slabs. In this example, we have $N=6$ and $M=4$.}
    \label{fig:intervals}
  \end{center}
\end{figure}

The following notation is used throughout this paper: Each component
$U_i(t)$, $i=1,\ldots,N$, of the approximate $\mathrm{m(c/d)G}(q)$
solution $U(t)$ of (\ref{eq:u'=f}) is a piecewise polynomial on a
partition of $(0,T]$ into $m_i$ sub-intervals.  Sub-interval $j$ for
component $i$ is denoted by $I_{ij}=(t_{i,j-1},t_{ij}]$, and the
length of the sub-interval is given by the local \emph{time step}
$k_{ij}=t_{ij}-t_{i,j-1}$.  We shall sometimes refer to $I_{ij}$ as an
\emph{element}.  This is illustrated in Figure~\ref{fig:intervals}.  On each
sub-interval $I_{ij}$, $U_{i}\vert_{I_{ij}}$ is a polynomial of degree
at most $q_{ij}$.

Furthermore, we shall assume that the interval $(0,T]$ is partitioned
into blocks between certain synchronized time levels
$0=T_0<T_1<\ldots<T_M=T$. For each $T_n$, $n = 0,1,\ldots,M$ and
each $i = 1,2,\ldots,N$, we require that there is a $0\leq j \leq
m_i$ such that $t_{ij} = T_n$. We refer to the collection of local
intervals between two synchronized time levels $T_{n-1}$ and $T_n$
as a \emph{time slab}.  We denote the length of a time slab by $K_n = T_n - T_{n-1}$.

\subsection{Outline of the paper}

We first give an introduction to multi-adaptive time-stepping in
Section~\ref{sec:intro}. We then present the key algorithms used by
the multi-adaptive ODE solver of DOLFIN in
Section~\ref{sec:algorithms}, followed by a discussion of data
structures for efficient representation and interpolation of
multi-adaptive solutions in Section~\ref{sec:datastructures}.  In
Section~\ref{sec:performance}, we discuss the efficiency of
multi-adaptive time-stepping and in Section~\ref{sec:examples}, we
present a number of numerical examples that demonstrate the efficiency
of the proposed algorithms and data structures. Finally, we give some
concluding remarks in Section~\ref{sec:conclusions}.

\section{Multi-adaptive time-stepping}
\label{sec:intro}

In this section, we give a quick introduction to multi-adaptive
time-stepping, including the formulation of the methods, error
estimates and adaptivity. For a more detailed account, we refer the
reader to~\cite{logg:article:01,logg:article:02,logg:article:08}.

\subsection{Formulation of the methods}

The \mcgq{} and \mdgq{} methods are obtained by multiplying the system of
equations~(\ref{eq:u'=f}) with a suitable test function~$v$, to obtain the
following variational problem:
Find $U\in V$ with $U(0)=u_0$, such that
\begin{equation} \label{eq:galerkin}
  \int_0^T (v, \dot{U}) \dt = \int_0^T (v, f(U,\cdot)) \dt \quad \forall v\in \hat{V},
\end{equation}
where $(\cdot,\cdot)$ denotes the standard inner product on $\real^N$
and $(\hat{V}, V)$ is a suitable pair of discrete function spaces,
the \emph{test} and \emph{trial} spaces respectively.

For the standard \cgq{} method, the trial space~$V$ consists of the
space of continuous piecewise polynomial vector-valued functions of
degree $q = q(t)$ on a partition $0 = t_0 < t_1 < \cdots < t_M = T$
and the test space~$\hat{V}$ consist of the space of (possibly
discontinuous) piecewise polynomial vector-valued functions of degree
$q - 1$ on the same partition. The multi-adaptive \mcgq{} method
extends the standard \cgq{} method by extending the test and trial
spaces to piecewise polynomial spaces on individual partitions of the
time interval that satisfy the constraints introduced in the previous
section and illustrated in Figure~\ref{fig:intervals}. Thus, each
component $U_i = U_i(t)$ is continuous and a piecewise polynomial on the
individual partition $0 = t_{i0} < t_{i1} < \cdots < t_{im_i} = T$ for
$i = 1,2,\ldots,N$.

For the standard \dgq{} method, the test and trial spaces are equal
and consist of the space of (possibly discontinuous) piecewise
polynomial vector-valued functions of degree $q = q(t)$ on a partition
$0 = t_0 < t_1 < \cdots < t_M = T$, which extends naturally to the
multi-adaptive \mdgq{} method by allowing each component of the test
and trial functions to be a piecewise polynomial on its own partition
of the time interval as above for the \mcgq{} method. Note that for
both the \dgq{} method and the \mdgq{} method, the integral
$\int_{0,T} (v, \dot{U}) \dt$ in~(\ref{eq:galerkin}) must be treated
appropriately at the points of discontinuity, see
\cite{logg:article:01}.

Both in the case of the \mcgq{} and \mdgq{} methods, the variational
problem~(\ref{eq:galerkin}) gives rise to a system of discrete
equations by expanding the solution~$U$ in a suitable basis on each
local interval~$I_{ij}$,
\begin{equation}
  U_i |_{I_{ij}} = \sum_{m=0}^{q_{ij}} \xi_{ijm} \phi_{ijm},
\end{equation}
where $\{\xi_{ijm}\}_{m=0}^{q_{ij}}$ are the \emph{degrees of freedom}
for $U_i$ on $I_{ij}$ and $\{\phi_{ijm}\}_{m=0}^{q_{ij}}$ is a
suitable basis for $P^{q_{ij}}(I_{ij})$. For any particular choice of
quadrature, the resulting system of discrete equations takes the form
of an implicit Runge--Kutta method on each local interval $I_{ij}$.
The discrete equations take the form
\begin{equation} \label{eq:equations}
  \xi_{ijm} =
  \xi_{ij0}^- +
  k_{ij} \sum_{n=0}^{q_{ij}} w_{mn}^{[q_{ij}]} \ f_i(U(\tau_{ij}^{-1}(s_{n}^{[q_{ij}]})),\tau_{ij}^{-1}(s_n^{[q_{ij}]})),
\end{equation}
for $m = 0,\ldots,q_{ij}$, where
$\{w_{mn}^{[q_{ij}]}\}_{m=0,n=0}^{q_{ij}}$ are weights, $\tau_{ij}$
maps $I_{ij}$ to $(0,1]$, $\tau_{ij}(t) =
(t-t_{i,j-1})/(t_{ij}-t_{i,j-1})$, and
$\{s_n^{[q_{ij}]}\}_{n=0}^{q_{ij}}$ are quadrature points defined on
$[0,1]$. Note that we have here assumed that the number of quadrature
points is equal to the number of nodal points.
See~\cite{logg:article:01} for a discussion of suitable quadrature
rules and basis functions.

\subsection{Error estimates and adaptivity}

The global error $e = U - u$ of the approximate solution $U$
of~(\ref{eq:u'=f}) may be bounded in terms of computable
quantities. Such an \emph{a posteriori} error estimate is proved
in~\cite{logg:article:01}, both for the \mcgq{} and \mdgq{} methods.
The a posteriori error estimate provides a bound for any given linear
functional $\mathcal{M} : \real^N \rightarrow \real$ of the global
error $e(T)$ at the final time, such as the error $e_i(T)$ in a single
component. Bounds for the error itself in various norms may also be
approximated. Below, we state the basic a posteriori error estimate
for the \mcgq{} method and refer to~\cite{logg:article:01} for a
complete discussion, including error estimates for~\mdgq{}.

For the \mcgq{} method, the error estimate takes the following form:
\begin{equation} \label{eq:estimate,cg}
  |\mathcal{M}(e(T))| \leq E \equiv \sum_{i=1}^N S_i(T) \max_{[0,T]} \left\{ C_i k_{i}^{q_{i}} |R_{i}| \right\},
\end{equation}
Here, $R = \dot{U} - f(U,\cdot)$ denotes the \emph{residual} of the
computed solution, $C_i = C_i(t)$ denotes an interpolation constant
(which may be different for each local interval) and $S_i(T)$ denotes
a \emph{stability factor} that measures the rate of propagation of
local errors for component $U_i$ (the influence of a nonzero residual
in component $U_i$ on the size of the error in the given
functional). By selecting the local time steps $k_i = k_i(t)$ such
that $E = \mathrm{TOL}$ for a given tolerance~$\mathrm{TOL}$, one may
thus guarantee that the error in the functional~$\mathcal{M}$ is
bounded by the given tolerance, $|\mathcal{M}(e(T))| \leq
\mathrm{TOL}$.

Comparing to standard Runge--Kutta methods for the
solution of initial value problems, the stability factor quantifies
the relationship between the ``local error'' and the global error. Note
that alternatively, the stability information may be kept as a local
time-dependent \emph{stability weight} for more fine-grained control
of the contributions to the global error. The stability factors are
obtained by solving a \emph{dual problem} of~(\ref{eq:u'=f}) for
the given functional~$\mathcal{M}$, see~\cite{EriEst95,logg:article:01}.
The particular form of the dual problem for~(\ref{eq:u'=f}) will be
discussed in Section~\ref{sec:dual}.

The individual time steps may be chosen so as to equidistribute
the error in the different components in an attempt to satisfy
\begin{equation} \label{eq:adaptivity}
  C_{ij} k_{ij}^{q_{ij}} \max_{I_{ij}} |R_i| = \mathrm{TOL} / (N S_i(T)),
\end{equation}
for each local time interval $I_{ij}$.
This may be done in an iterative fashion, as outlined in the following
basic adaptive algorithm:
\begin{enumerate}
\item[(0)\,]
  Assume $S_i(T) = 1$ for $i = 1,2,\ldots,N$;
\item[(i)\,\,\,]
  Solve the primal problem with time steps based on (\ref{eq:adaptivity});
\item[(ii)\,]
  Solve the dual problem and compute the stability factors;
\item[(iii)]
  Compute an error bound $E$ based on (\ref{eq:estimate,cg});
\item[(iv)]
  If $E \leq \mathrm{TOL}$ then stop; if not go back to (i).
\end{enumerate}

\section{Algorithms}
\label{sec:algorithms}

We present below a collection of key algorithms for multi-adaptive
time-stepping. The algorithms are given in pseudo-code and where
appropriate we give remarks on how the algorithms have been
implemented in C++ for DOLFIN. In most cases, we present simplified
versions of the algorithms with focus on the most essential steps.

\subsection{General algorithm}

The general multi-adaptive time-stepping algorithm is Algorithm
\ref{alg:integrate}. Starting at $t = 0$, the algorithm creates a
sequence of time slabs until the given end time $T$ is reached. In
each macro time step, Algorithm~\ref{alg:createtimeslab}
(CreateTimeSlab) is called to create a time slab covering an interval
$[T_{n-1},T_n]$ such that $T_n \leq T$. For each time slab, the system
of discrete equations is solved iteratively, using direct fixed-point
iteration or a preconditioned Newton's method, until the discrete
equations given by the \mcgq{} or \mdgq{} method have converged.

\begin{algorithm}
  \begin{tabbing}
    $t \leftarrow 0$ \\
    \textbf{while} {$t < T$} \\
    \tab \{time slab, $t$\} $\leftarrow$ CreateTimeSlab($\{1,\ldots,N\}$, $t$, $T$) \\
    \tab SolveTimeSlab(time slab) \\
    \textbf{end while}
  \end{tabbing}
  \caption{$U = $ Integrate(ODE)}
  \label{alg:integrate}
\end{algorithm}

The basic forward integrator, Algorithm \ref{alg:integrate}, can be
used as the main component of an adaptive algorithm with automated
error control of the computed solution as outlined in
Section~\ref{sec:intro}.  In each iteration, the \emph{primal problem}
(\ref{eq:u'=f}) is solved using Algorithm \ref{alg:integrate}. An ODE
of the form (\ref{eq:u'=f}) representing the \emph{dual problem} is
then created and solved using Algorithm \ref{alg:integrate}.  It is
important to note that both the primal and the dual problems may be
solved using the same algorithm, but with (possibly) different time
steps, tolerances, methods, and orders.  When the solution of the dual
problem has been computed, the stability factors $\{S_i(T)\}_{i=1}^N$
and the error estimate may be computed.

\subsection{Recursive construction of time slabs}
\label{sec:construction}

\begin{figure}[htbp]
  \begin{center}
    \psfrag{t1}{$T_{n-1}$}
    \psfrag{t2}{$T_n$}
    \includegraphics[width=12cm]{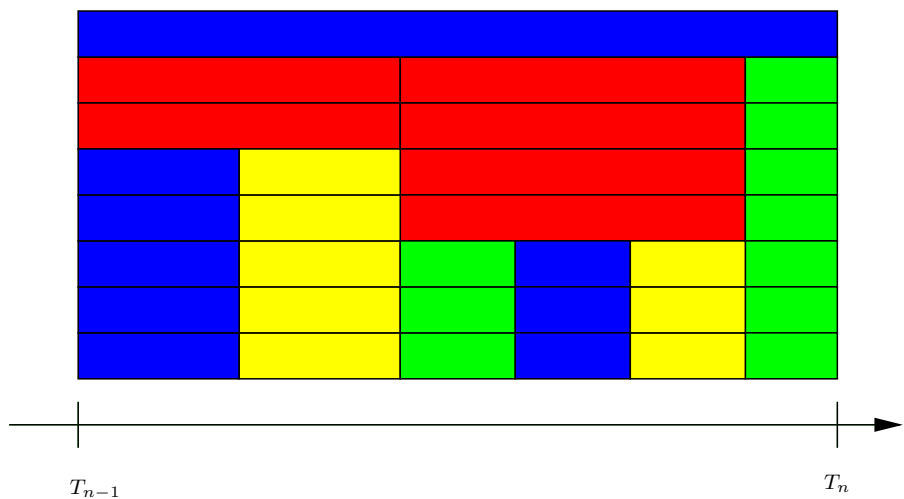}
    \caption{The recursive organization of the time slab. Each time
      slab contains an element group and a list of recursively nested
      time slabs. The root time slab in the figure contains one
      element group of one element and three sub-slabs. The first
      of these sub-slabs contains an element group of two elements and
      two nested sub-slabs, and so on. The root time slab recursively
      contains a total of nine element groups and 33 elements.}
    \label{fig:recursiveslab}
  \end{center}
\end{figure}

In each step of Algorithm \ref{alg:integrate}, a new time slab is
created between two synchronized time levels $T_{n-1}$ and $T_n$.  The
time slab is organized recursively as follows. The root time slab
covering the interval $[T_{n-1}, T_n]$ contains a non-empty list of
elements, which we refer to as an \emph{element group}, and a possibly
empty list of time slabs, which in turn may contain nested groups of
elements and time slabs. Each such element group together with the
corresponding nested set of element groups is referred to as a
\emph{sub-slab}. This is illustrated in Figure \ref{fig:recursiveslab}.

\begin{figure}[htbp]
  \begin{center}
    \psfrag{k1}{$\theta K$}
    \psfrag{k2}{$\underline{K}$}
    \psfrag{k3}{$K$}
    \includegraphics[width=10cm]{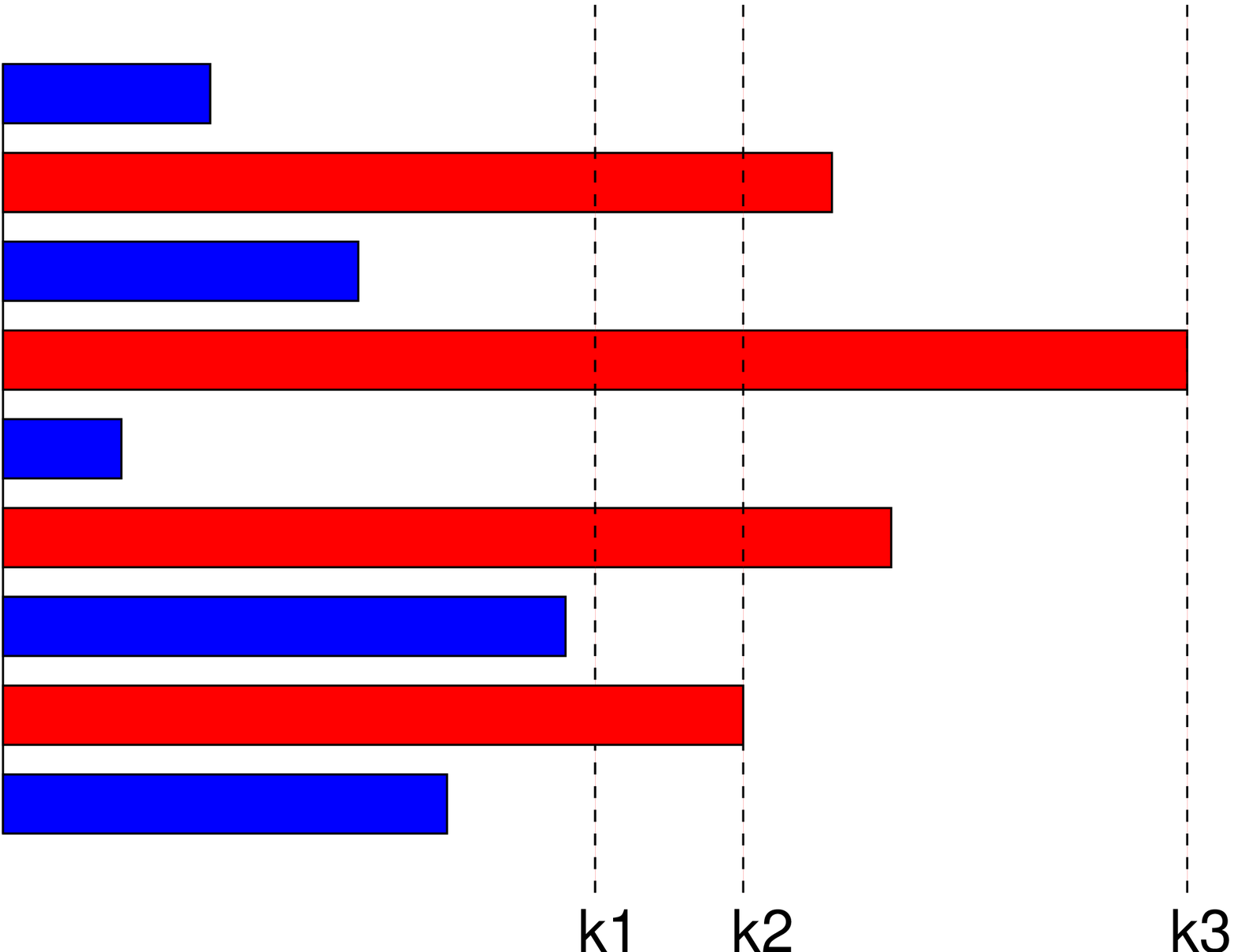}
    \caption{The partition of components into groups of small and
      large time steps for $\theta = 1/2$.}
    \label{fig:partition}
  \end{center}
\end{figure}

To create a time slab, we first compute the desired time steps for all
components as given by the a posteriori error
estimate~(\ref{eq:estimate,cg}). We discuss in detail the time step
selection below in Section~\ref{sec:adaptivity}. A threshold $\theta
K$ is then computed based on the maximum time step $K$ among the
components and a fixed parameter $\theta\in(0,1)$ controlling the
density of the time slab. The components are partitioned into two sets
based on the threshold, and a large time step $\underline{K}$ is
selected to be the smallest time step among the components in the set
with large time steps as described in Algorithm~\ref{alg:partition}
and illustrated in Figure \ref{fig:partition}. For each component in
the group with large time steps, an element is created and added to
the element group of the time slab.  The remaining components with
small time steps are processed by a recursive application of this
algorithm for the construction of time slabs.

We organize the recursive construction of time slabs as described by
Algorithms \ref{alg:createtimeslab}, \ref{alg:partition},
\ref{alg:createelements}, and \ref{alg:createtimeslabs}.  The
recursive construction simplifies the implementation; each recursively
nested sub-slab can be considered as a sub-system of the ODE. Note
that the element group containing elements for components in group
$I_1$ is created before the recursively nested sub-slabs for
components in group $I_0$. The tree of time slabs is thus created
recursively \emph{breadth-first}, which means in particular that the
element for the component with the largest time step is created first.

Algorithm \ref{alg:partition} for the partition of components can be
implemented efficiently using the function \texttt{std::partition()}, which
is part of the Standard C++ Library.

\begin{algorithm}
  \begin{tabbing}
    \{$I_0$, $I_1$, $K$\} $\leftarrow$ Partition(components) \\
    \textbf{if} {$T_{n-1} + K < T$} \\
    \tab $T_n \leftarrow T_{n-1} + K$ \\
    \textbf{else} \\
    \tab $T_n \leftarrow T$ \\
    \textbf{end if} \\
    element group $\leftarrow$ CreateElements($I_1$, $T_{n-1}$, $T_n$) \\
    time slabs $\leftarrow$ CreateTimeSlabs($I_0$, $T_{n-1}$, $T_n$) \\
    time slab $\leftarrow$ \{element group, time slabs\}
  \end{tabbing}
  \caption{\{time slab, $T_n$\} = CreateTimeSlab(components, $T_{n-1}$, $T$)}
  \label{alg:createtimeslab}
\end{algorithm}

\begin{algorithm}
  \begin{tabbing}
    $I_0 \leftarrow \emptyset$ \\
    $I_1 \leftarrow \emptyset$ \\
    $K \leftarrow$ \emph{maximum time step within} components \\
    \textbf{for each} component \\
    \tab $k \leftarrow$ \emph{time step of} component \\
    \tab \textbf{if} {$k < \theta K$} \\
    \tab \tab $I_0 \leftarrow I_0$ $\cup$ \{component\} \\
    \tab \textbf{else} \\
    \tab \tab $I_1 \leftarrow I_1$ $\cup$ \{component\} \\
    \tab \textbf{endif} \\
    \textbf{end for} \\
    $\underline{K} \leftarrow$ \emph{minimum time step within} $I_1$ \\
    $K \leftarrow \underline{K}$
  \end{tabbing}
  \caption{\{$I_0$, $I_1$, $K$\} = Partition(components)}
  \label{alg:partition}
\end{algorithm}

\begin{algorithm}
  \begin{tabbing}
     elements $\leftarrow \emptyset$ \\
     \textbf{for} {\textbf{each} component} \\
     \tab \emph{create} element \emph{for} component \emph{on} $[T_{n-1},T_n]$ \\
     \tab elements $\leftarrow$ elements $\cup$ element \\
     \textbf{end for}
  \end{tabbing}
  \caption{elements = CreateElements(components, $T_{n-1}$, $T_n$)}
  \label{alg:createelements}
\end{algorithm}

\begin{algorithm}
  \begin{tabbing}
     time slabs $\leftarrow \emptyset$ \\
     $t \leftarrow T_{n-1}$ \\
     \textbf{while} {$t < T$} \\
     \tab \{time slab, $t$\} $\leftarrow$ CreateTimeSlab(components, $t$, $T_n$) \\
     \tab time slabs $\leftarrow$ time slabs $\cup$ time slab \\
     \textbf{end while}
  \end{tabbing}
  \caption{time slabs = CreateTimeSlabs(components, $T_{n-1}$, $T_n$)}
  \label{alg:createtimeslabs}
\end{algorithm}

\subsection{Solving the system of discrete equations}
\label{sec:iteration}

On each time slab $\mathcal{T}_n$, $n = 1,2,\ldots,M$, we need to
solve a system of equations for the degrees of freedom on the time
slab. On each local interval $I_{ij} \in \mathcal{T}_n$, these equations are
given by~(\ref{eq:equations}). Depending on the properties of the
given system~(\ref{eq:u'=f}), different solution strategies for the
time slab system~(\ref{eq:equations}) may be appropriate as outlined
below.

\subsubsection{Direct fixed-point iteration}

In the simplest case, the time slab system is solved by direct
fixed-point iteration on~(\ref{eq:equations}) for each element in the
time slab. The fixed-point iteration is performed in a forward
fashion, sweeping over the elements in the time slab in the same order
as they are created by Algorithm~\ref{alg:createtimeslab}. In
particular, this means that for each component in the time slab
system, the end-time value on each element is updated before the
degrees of freedom for the following element. Thus, for each element
$I_{ij} \in \mathcal{T}_n$, we compute the degrees of freedom
$\{\xi_{ijm}\}_{j=0}^{q_{ij}}$ according to
\begin{equation} \label{eq:equations,again}
  \xi_{ijm} =
  \xi_{ij0}^- +
  k_{ij} \sum_{n=0}^{q_{ij}} w_{mn}^{[q_{ij}]} \
  f_i(U(\tau_{ij}^{-1}(s_{n}^{[q_{ij}]})),\tau_{ij}^{-1}(s_n^{[q_{ij}]})),
  \quad m = 0,1,\ldots,q_{ij}.
\end{equation}
Direct fixed-point iteration
converges if the system is non-stiff and typically only a few
iterations are needed. In fact, one may consider a system to be stiff
if direct fixed-point iteration does not converge.

\subsubsection{Damped fixed-point iteration}

If the system is stiff, that is, direct fixed-point iteration does not
converge, one may introduce a suitable amount of damping to adaptively
stabilize the fixed-point iteration. The fixed-point
iteration~(\ref{eq:equations,again}) may be written in the form
\begin{equation} \label{eq:fixed-point}
  \xi_{ijm} = g_{ijm}(\xi),
\end{equation}
where $\xi$ is the vector of degrees of freedom for the
solution on the time slab. We modify the fixed-point iteration by
introducing a damping parameter $\alpha$:
\begin{equation} \label{eq:fixed-point,damped}
  \xi_{ijm} = (1-\alpha_{ijm}) \xi_{ijm} + \alpha_{ijm} g_{ijm}(\xi).
\end{equation}
In \cite{logg:thesis:03}, a number of different strategies for the
selection of the damping parameter $\alpha$ are discussed.  We mention
two of these strategies here. The first strategy chooses $\alpha$
based on the diagonal derivatives $\partial f_i/\partial u_i$, $i =
1,2,\ldots,N$, corresponding to a modified Newton's method where the
Jacobian is approximated by a diagonal matrix. This strategy works
well for systems with a diagonally dominant Jacobian, including many
systems arising when modeling chemical reactions. The second strategy
adaptively chooses a scalar $\alpha$ based on the convergence of the
fixed-point iterations.

\subsubsection{Newton's method}

Alternatively, one may apply Newton's method directly to the full
system of equations~(\ref{eq:equations,again})
associated with each time slab. The linear system in each Newton
iteration may then be solved either by a direct method or an iterative
method such as a Krylov subspace method in combination with a suitable
preconditioner, depending on the characteristics of the underlying
system~(\ref{eq:u'=f}). In addition, one may also apply a special
preconditioner that improves the convergence by propagating values
forward in time within the time slab. Note that if the multi-adaptive
efficiency index is large (see Section~\ref{sec:performance} below),
then the time slab system is not significantly larger than the
corresponding time slab system for a mono-adaptive method.

\subsubsection{Choosing a solution strategy}

Ultimately, an intelligent solver should automatically choose a
suitable algorithm for the solution of the time slab system. Thus, the
solver may initially try direct fixed-point iteration. If the system
is stiff, the solver switches to adaptive fixed-point iteration (as
outlined in~\cite{logg:thesis:03}). Finally, if the adaptive fixed-point
iteration converges slowly, the solver may switch to Newton's method.

\subsubsection{Interpolation of the solution}

To update the degrees of freedom on an element according to
(\ref{eq:equations,again}), the appropriate component $f_i$ of the
right-hand side of (\ref{eq:u'=f}) needs to be evaluated at the set of
quadrature points. In order for $f_i$ to be evaluated, each component
$U_{i'}$ of the computed solution $U$ on which $f_i$ depends, needs to be
evaluated at the quadrature points. We let $\mathcal{S}_i \subseteq
\{1,\ldots,N\}$ denote the \emph{sparsity pattern} of component $U_i$,
that is, the set of components on which $f_i$ depends,
\begin{equation}
  \mathcal{S}_i = \{ i' \in \{1,\ldots,N\} :
  \partial f_i / \partial u_{i'} \neq 0 \}.
\end{equation}
Thus, to evaluate $f_i$ at a given quadrature point $t$, only the
components $\{U_{i'}\}_{i' \in \mathcal{S}_i}$ need to be evaluated at $t$, as
in Algorithm \ref{alg:rhs}. This is of particular importance for
problems of sparse structure and enables efficient multi-adaptive
integration of time-dependent PDEs, as demonstrated below in
Section~\ref{sec:examples}.  The sparsity pattern $\mathcal{S}_i$ is
automatically detected by the solver. Alternatively, the sparsity
pattern may be specified by a (sparse) matrix.

\begin{algorithm}
  \begin{tabbing}
    \textbf{for} {$i' \in \mathcal{S}_i$} \\
    \tab $x(i') \leftarrow U_{i'}(t)$ \\
    \textbf{end for} \\
    $y \leftarrow f_i(x,t)$
  \end{tabbing}
  \caption{$y = $ EvaluateRightHandSide($i$, $t$)}
  \label{alg:rhs}
\end{algorithm}

In Algorithm~\ref{alg:rhs}, the key step is the evaluation of a
component $U_{i'}$ at a given point $t$. For a standard mono-adaptive
method, this is straightforward since all components use the same time
steps. In particular, if the quadrature points are chosen to be the
same as the nodal points, the value of $U_{i'}(t)$ is
known. For a multi-adaptive method, a quadrature point~$t$ for the
evaluation of $f_i$ is not necessarily a nodal point for $U_{i'}$. To
evaluate $U_{i'}(t)$, one thus needs to find the local interval
$I_{i'j'}$ such that $t\in I_{i'j'}$ and then evaluate $U_{i'}(t)$ by
interpolation on that interval. In Section~\ref{sec:datastructures}
below, we discuss data structures that allow efficient storage and
interpolation of the multi-adaptive solution. In particular, these
data structures give $\mathcal{O}(1)$ access to the value of any
component $U_{i'}$ in the sparsity pattern $\mathcal{S}_i$ at any
quadrature point $t$ for $f_i$.

\subsection{Multi-adaptive time step selection}
\label{sec:adaptivity}

The individual and adaptive time steps $k_{ij}$ are determined during
the recursive construction of time slabs based on an a posteriori
error estimate as discussed in Section~\ref{sec:intro}. Thus,
according to~(\ref{eq:adaptivity}), each local time step $k_{ij}$
should be chosen to satisfy
\begin{equation} \label{eq:timesteps}
  k_{ij} =
  \left(
  \frac{\mathrm{TOL}}{C_{ij} N S_i(T) \max_{I_{ij}} |R_i|}
  \right)^{1/q_{ij}}.
\end{equation}
where $\mathrm{TOL}$ is a given tolerance.

However, the time steps can not be based directly on
(\ref{eq:timesteps}), since that leads to unwanted oscillations in the
size of the time steps. If $r_{i,j-1} = \max_{I_{i,j-1}} |R_i|$ is
small, then $k_{ij}$ will be large, and as a result $r_{ij}$ will also
be large. Consequently, $k_{i,j+1}$ and $r_{i,j+1}$ will be small, and
so on. To avoid these oscillations, we adjust the time step $k_{ij}$
according to Algorithm \ref{alg:regulator}, which determines the new
time step as a weighted harmonic mean value of the previous time step
and the time step given by (\ref{eq:timesteps}).  Alternatively,
DOLFIN provides time step control based on the PID controllers
presented in~\cite{GusLun88,Sod03}, including H0211 and
H211PI. However, the simple controller of
Algorithm~\ref{alg:regulator} performs well compared to the more
sophisticated controllers in~\cite{GusLun88,Sod03}. A
suitable value for the weight $w$ in Algorithm~\ref{alg:regulator} is
$w = 5$ (found empirically).

\begin{algorithm}
  \begin{tabbing}
    $k \leftarrow (1 + w) k_{\mathrm{old}} k_{\mathrm{new}} / (k_{\mathrm{old}} + wk_{\mathrm{new}})$ \\
    $k \leftarrow \min(k, k_{\max})$
  \end{tabbing}
  \caption{$k = $ Controller($k_{\mathrm{new}}$, $k_{\mathrm{old}}$, $k_{\max}$)}
  \label{alg:regulator}
\end{algorithm}

The initial time steps $k_{11} = k_{21} = \cdots = k_{N1} = K_1$ are
chosen equal for all components and are determined iteratively for the
first time slab. The size $K_1$ of the first time slab is first
initialized to some default value, possibly based on the length $T$ of
the time interval, and then adjusted until the local residuals are
sufficiently small for all components.

\subsection{Solving the dual problem}
\label{sec:dual}

Stability factors may be approximated by numerically solving an
auxiliary dual problem for (\ref{eq:u'=f}). This dual problem is
given by the following system of linear ordinary differential equations:
\begin{equation}
  \label{eq:dual}
  \begin{split}
    - \dot{\varphi}(t) &= J(U(t),t)^{\top} \varphi(t), \quad t\in [0,T), \\
      \varphi(T) &= \psi,
  \end{split}
\end{equation}
where $J(U(t),t)$ denotes the Jacobian of the right-hand side $f$ of
(\ref{eq:u'=f}) at time $t$ and $\psi = \mathcal{M}'$ (the Riesz
representer of $\mathcal{M}$) is initial data for the dual problem
corresponding to the given functional $\mathcal{M}$ to be estimated.
Note that we need to linearize around the computed solution $U$, since
the exact solution $u$ of (\ref{eq:u'=f}) is not known.  To solve this
backward problem over $[0,T)$ using the forward integrator Algorithm
  \ref{alg:integrate}, we rewrite (\ref{eq:dual}) as a forward
  problem. With $w(t) = \varphi(T-t)$, we have $\dot{w} =
  -\dot{\varphi}(T-t) = J(U(T-t),T-t)^{\top} w(t)$, and so (\ref{eq:dual})
  can be written as a forward problem for $w$ in the form
\begin{equation}
  \label{eq:dual,w}
  \begin{split}
    \dot{w}(t) &= f^*(w(t),t) \equiv J(U(T-t),T-t)^{\top} w(t), \quad t\in (0,T], \\
    w(0) &= \psi.
  \end{split}
\end{equation}

\section{Data structures}
\label{sec:datastructures}

For a standard mono-adaptive method, the solution on a time slab is
typically stored as an array of values at the right end-point of the
time slab, or as a list of arrays (possibly stored as one contiguous
array) for a higher order method with several stages. However, a
different data structure is needed to store the solution on a
multi-adaptive time slab. Such a data structure should ideally store
the solution with minimal overhead compared to the cost of storing
only the array of degrees of freedom for the solution on the time
slab. In addition, it should also allow for efficient interpolation of
the solution, that is, accessing the values of the solution for all
components at any given time within the time slab. We present below a
data structure that allows efficient storage of the entire solution on
a time slab with little overhead, and at the same time allows
efficient interpolation with $\mathcal{O}(1)$ access to any given
value during the iterative solution of the system of discrete
equations.

\subsection{Representing the solution}

The multi-adaptive solution on a time-slab can be efficiently
represented using a data structure consisting of eight arrays as
shown in Table~\ref{tab:datastructure}. For simplicity, we assume that
all elements in a time slab are constructed for the same choice of
method, \mcgq{} or \mdgq{}, for a given fixed $q$.

The recursive construction of time slabs as discussed in
Section~\ref{sec:construction} generates a sequence of \emph{sub
slabs}, each containing a list of \emph{elements} (an element group).
For each sub-slab, we store the value of the time $t$ at the left
end-point and at the right end-point in the two arrays \texttt{sa} and
\texttt{sb}. Thus, for sub-slab number~$s$ covering the interval~$(a_s,b_s)$, we have
\begin{equation}
  \begin{split}
    a_s &= \texttt{sa[}s\texttt{]}, \\
    b_s &= \texttt{sb[}s\texttt{]}.
  \end{split}
\end{equation}
Furthermore, for all elements in the (root) time slab, we store the
degrees of freedom in the order they are created in the array
\texttt{jx} (mapping a degree of freedom $j$ to the value $x$ of that
degree of freedom). Thus, if each element has $q$ degrees of freedom,
as in the case of the multi-adaptive $\mathrm{mcG}(q)$ method, then
the length of the array \texttt{jx} is $q$ times the number of
elements. In particular, if all components use the same time steps,
then the length of the array \texttt{jx} is $qN$.

For each element, we store the corresponding component index $i$ in
the array \texttt{ei} in order to be able to evaluate the correct
component $f_i$ of the right-hand side $f$ of (\ref{eq:u'=f}) when
iterating over all elements in the time slab to update the degrees of
freedom. When updating the values on an element according to
(\ref{eq:equations,again}), it is also necessary to know the left and right
end-points of the elements. Thus, we store an array \texttt{es} that
maps the number $e$ of a given element to the number $s$ of the
corresponding sub-slab containing the element. As a consequence, the
left end-point $a_e$ and right end-point $b_e$ for a given element $e$
are given by
\begin{equation}
  \begin{split}
    a_e &= \mathtt{sa[es[}e\mathtt{]]}, \\
    b_e &= \mathtt{sb[es[}e\mathtt{]]}.
  \end{split}
\end{equation}

\begin{table}[htbp]
  \begin{center}
    \begin{tabular}{lll}
      \hline
      Array & Type & Description \\
      \hline
      \hline
      \texttt{sa} & \texttt{double} & left end-points for sub-slabs \\
      \texttt{sb} & \texttt{double} & right end-points for sub-slabs \\
      \texttt{jx} & \texttt{double} & values for degrees of freedom \\
      \texttt{ei} & \texttt{int}    & component indices for elements \\
      \texttt{es} & \texttt{int}    & time slabs containing elements \\
      \texttt{ee} & \texttt{int}    & previous elements for elements \\
      \texttt{ed} & \texttt{int}    & first dependencies for elements \\
      \texttt{de} & \texttt{int}    & elements for dependencies \\
      \hline
    \end{tabular}
    \caption{Data structures for efficient representation of a
      multi-adaptive time slab.}
    \label{tab:datastructure}
  \end{center}
\end{table}

\subsection{Interpolating the solution at quadrature points}

Updating the values on an element according to (\ref{eq:equations,again})
also requires knowledge of the value at the left end-point, which is
given as the end-time value on the previous element in the time slab
for the same component (or the end-time value from the previous time
slab). This information is available in the array \texttt{ee}, which
stores for each element the number of the previous element (or $-1$ if
there is no previous element).

As discussed above in Section~\ref{sec:iteration}, the system of
discrete equations on each time slab is solved by iterating over the
elements in the time slab and updating the values on each element,
either in a direct fixed-point iteration or a Newton's method. We must
then for any given element~$e$ corresponding to some component $i =
\mathtt{ei[}e\mathtt{]}$ evaluate the right-hand side $f_i$ at each
quadrature point $t$ within the element. This requires the values of
the solution $U$ at $t$ for all components contained in the
sparsity pattern $\mathcal{S}_i$ for component $i$ according to
Algorithm~\ref{alg:rhs}. As a consequence of
Algorithm~\ref{alg:createtimeslab} for the recursive construction of
time slabs, elements for components that use large time steps are
constructed before elements for components that use small time
steps. Since all elements of the time slab are traversed in the same
order during the iterative solution of the system of discrete
equations, elements corresponding to large time steps have recently
been visited and cover any element that corresponds to a smaller time
step. The last visited element for each component is stored in an
auxiliary array~\texttt{elast} of size~$N$. Thus, if $i' \in \mathcal{S}_i$ and
component $i'$ has recently been visited, then it is straight-forward
to find the latest element $e' = \mathtt{elast}[i']$ for component $i'$
that covers the current element for component $i$ and interpolate
$U_{i'}$ at time $t$. It is also straight-forward to interpolate the
values for any components that are present in the same element group
as the current element.

However, when updating the values on an element $e$ corresponding to
some component $i = \mathtt{ei[}e\mathtt{]}$ depending on some other
component $i' \in \mathcal{S}_i$ which uses smaller time steps, one
must find for each quadrature point $t$ on the element $e$ the element
$e'$ for component $i'$ containing $t$, which is non-trivial. The
element $e'$ can be found by searching through all elements for
component $i'$ in the time slab, but this quickly becomes
inefficient. Instead, we store for each element $e$ a list of
dependencies to elements with smaller time steps in the two arrays
\texttt{ed} and \texttt{de}. These two arrays store a sparse integer
matrix of dependencies to elements with smaller time steps for all
elements in the time slab. Thus, for any given element $e$, the number
of dependencies to elements with smaller time steps is given by
\begin{equation}
  \mathtt{ed[}e + 1\mathtt{]} - \mathtt{ed[}e\mathtt{]},
\end{equation}
and the elements with smaller time steps that need to be interpolated
at the quadrature points for element~$e$ are given by
\begin{equation}
  \{
  \mathtt{de[ed[}e\mathtt{]]},
  \mathtt{de[ed[}e\mathtt{]} + 1\mathtt{]},
  \ldots,
  \mathtt{de[ed[}e+1\mathtt{]} -1 \mathtt{]}
  \}.
\end{equation}

\section{Performance}
\label{sec:performance}

The efficiency of multi-adaptive time-stepping compared to standard
mono-adaptive time-stepping depends on the system being integrated,
the tolerance, and the efficiency of the implementation. For many
systems, the potential speedup is large, but the actual speedup
depends also on the overhead needed to handle the additional
complications of a multi-adaptive implementation: the recursive
construction of time slabs and the interpolation of values within a
time slab.

To study the performance of multi-adaptive time-stepping, we consider
a system of $N$ components and time steps given by $\{k_{ij} =
|I_{ij}| : I_{ij} \in \mathcal{T}_n\}$ on some time slab
$\mathcal{T}_n$. We define the \emph{multi-adaptive efficiency index}~$\mu$ by
\begin{equation}
\mu =
\frac{N / k_{\min}}{|\mathcal{T}_n| / k_{\max}}
=
\frac{k_{\max}}{k_{\min}} \,
\frac{N}{|\mathcal{T}_n|},
\end{equation}
where $k_{\min} = \min_{I_{ij}\in\mathcal{T}_n} k_{ij}$,
$k_{\max} = \max_{I_{ij}\in\mathcal{T}_n} k_{ij}$ and
$|\mathcal{T}_n|$ is the number of local intervals in the time slab $\mathcal{T}_n$.
Thus, to obtain the multi-adaptive efficiency index,
we divide the number of local intervals per unit time for a
mono-adaptive discretization with the actual number of local intervals
per unit time for a multi-adaptive discretization. This is the
potential speedup when compared to a mono-adaptive method that is
forced to use the same small time step $k_{\min}$ for all
components. However, the actual speedup is always smaller than~$\mu$
for two reasons. The first is the overhead of the multi-adaptive
implementation and the second is that the system of discrete equations
on each time slab may sometimes be more expensive to solve than the
corresponding mono-adaptive systems (because they are typically larger
in size).

\begin{figure}[htbp]
  \begin{center}
    \psfrag{t1}{$T_{n-1}$}
    \psfrag{t2}{$T_n$}
    \includegraphics[width=12cm]{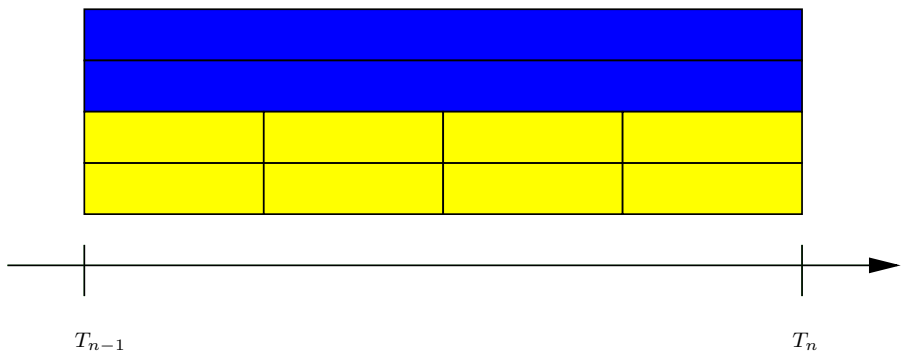}
    \caption{A time slab with $N_K = N_k = 2$ and multi-adaptive efficiency index~$\mu = 16/10=1.6$}.
    \label{fig:model}
  \end{center}
\end{figure}

Consider a model problem consisting of $N = N_K + N_k$ components, where $N_K$
components vary on a slow time scale $K$ and $N_k$ components vary on
a fast time scale $k$ as in Figure~\ref{fig:model}. The potential speedup is given by the
multi-adaptive efficiency index,
\begin{equation} \label{eq:efficiencyindex}
  \mu =
  \frac{K}{k} \, \frac{N}{N_K + N_k K/k}
  =
  \frac{K}{k} \, \frac{N/K}{N_K/K + N_k/k}
  \sim
  \frac{K}{k} \gg 1,
\end{equation}
if $N_K/K \gg N_k/k$ and $K \gg k$, that is the number of large elements dominates
the number of small elements. Thus, the potential speedup can be very
large for a system where a large part of the system varies on a large
time scale and a small part of the system varies on a small time
scale.

If, on the other hand, $K \sim k$ or $N_K \sim N_k$, then the
multi-adaptive efficiency index may be of moderate size. As a
consequence, the actual speedup may be small (or even ``negative'') if
the overhead of the multi-adaptive implementation is significant. In
the next section, we indicate the multi-adaptive efficiency index and
compare this to the actual speedup for a number of benchmark problems.

\section{Numerical examples and benchmark results}
\label{sec:examples}

In this section, we present two benchmark problems to demonstrate
the efficiency of multi-adaptive time-stepping. Both examples are
time-dependent PDEs that we discretize in space using the
$\mathrm{cG}(1)$ finite element method to obtain a system of ODEs,
sometimes referred to as the method of lines approach. In each case,
we lump and invert the mass matrix so as to obtain a system of the
form~(\ref{eq:u'=f}).

In the first of the two benchmark problems, the individual time steps
are chosen automatically based on an a~posteriori error estimate as
discussed above in Section~\ref{sec:adaptivity}. For the second
problem, the time steps are fixed in time and determined according to
a local CFL condition $k \sim h$ on each element. The results were
obtained with DOLFIN version 0.6.2.

\subsection{A nonlinear reaction-diffusion equation}

As a first example, we solve the following nonlinear
reaction-diffusion equation, taken from~\cite{SavHun05}:
\begin{equation} \label{eq:reaction}
  \begin{split}
    u_t - \epsilon u_{xx} &= \gamma u^2 (1 - u) \quad \mbox{ in } \Omega \times (0,T], \\
    \partial_n u &= 0 \quad \mbox{ on } \partial \Omega \times (0,T], \\
    u(\cdot, 0) &= u_0 \quad \mbox{ in } \Omega,
  \end{split}
\end{equation}
with $\Omega = (0,L)$, $\epsilon = 0.01$, $\gamma = 1000$ and final time $T = 1$.

The equation is discretized in space with the standard
$\mathrm{cG}(1)$ method using a uniform mesh with $1000$~mesh points.
The initial data is chosen according to
\begin{equation}
  u_0(x) = \frac{1}{1 + \exp(\lambda(x - 1))}.
\end{equation}
The resulting solution is a reaction front, sweeping across the domain
from left to right, as demonstrated in Figure~\ref{fig:reaction,u}.
The multi-adaptive time steps are automatically selected to be small
in and around the reaction front and sweep the domain at the same
velocity as the reaction front, as demonstrated in
Figure~\ref{fig:reaction,k}.

\begin{figure}[htbp]
  \begin{center}
    \psfrag{x}{$x$}
    \psfrag{u}{$U(x)$}
    \psfrag{t1}{\hspace{-0.2cm}\scriptsize $t = 0$}
    \psfrag{t2}{\hspace{-0.2cm}\scriptsize $t = 0.25$}
    \psfrag{t3}{\hspace{-0.2cm}\scriptsize $t = 0.5$}
    \psfrag{t4}{\hspace{-0.2cm}\scriptsize $t = 0.75$}
    \psfrag{t5}{\hspace{-0.2cm}\scriptsize $t = 1$}
    \includegraphics[width=12cm]{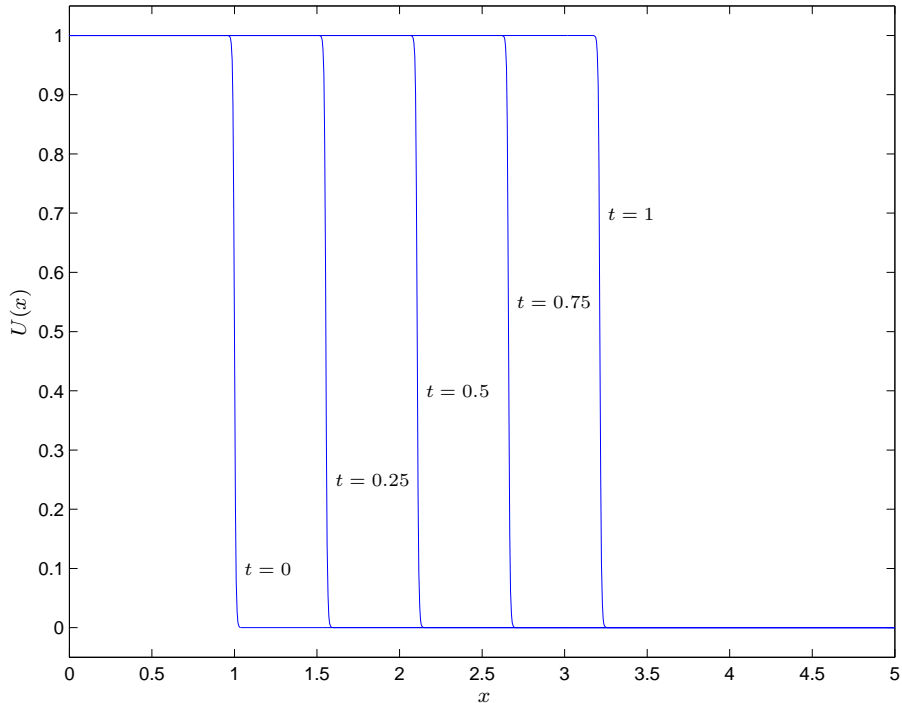}
    \caption{Propagation of the solution of the reaction--diffusion problem~(\ref{eq:reaction}).}
    \label{fig:reaction,u}
  \end{center}
\end{figure}

\begin{figure}[htbp]
  \begin{center}
    \psfrag{x}{$x$}
    \psfrag{t}{$t$}
    \psfrag{k}{$k(x,t)$}
    \includegraphics[width=10cm]{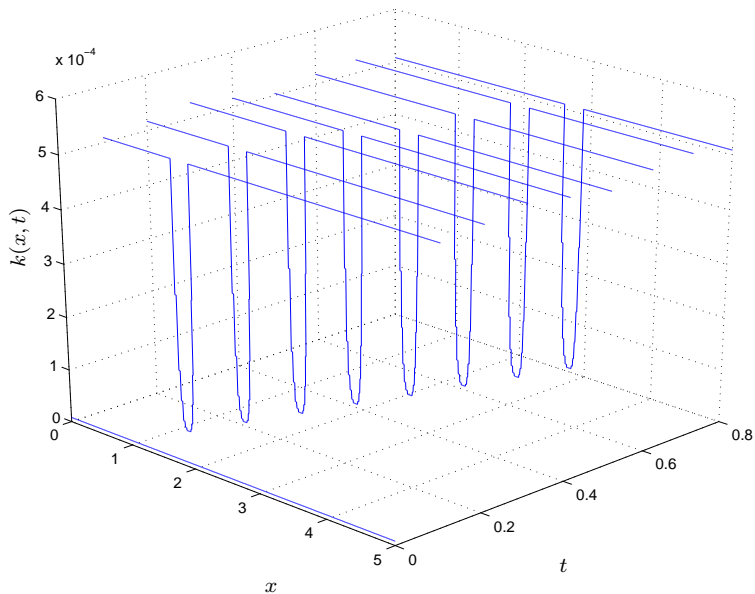}
    \caption{The multi-adaptive time steps as function of space at a
      sequence of points in time for the test problem~(\ref{eq:reaction}).}
    \label{fig:reaction,k}
  \end{center}
\end{figure}

To study the performance of the multi-adaptive solver, we compute the
solution for a range of tolerances with $L = 5$ and compare the
resulting error and CPU time with a standard mono-adaptive solver that
uses equal (adaptive) time steps for all components. To make the
comparison fair, we compare the multi-adaptive $\mathrm{mcG}(q)$
method with the mono-adaptive $\mathrm{cG}(q)$ method. In the
benchmarks, we only examine $q=1$. Both methods are implemented for
general order $q$ in the same programming language (C++) within a
common framework (DOLFIN), but the mono-adaptive method takes full
advantage of the fact that the time steps are equal for all
components. In particular, the mono-adaptive solver may use much
simpler data structures (a plain C array) to store the solution on
each time slab and there is no overhead for interpolation of the
solution. Furthermore, for the multi-adaptive solver, we need to
supply a right-hand side function $f$ which may be called to evaluate
single components $f_i(U(t), t)$, while for the mono-adaptive solver,
we may evaluate all components of $f$ at the same time, which is
usually an advantage (for the mono-adaptive solver).

This is a more meaningful measure of performance compared to only
measuring the number degrees of freedom (local steps) or comparing the
CPU time against the same multi-adaptive solver when it is forced to
use identical time steps for all components as
in~\cite{logg:article:01}, since one must also take into account the
overhead of the more complicated algorithms and data structures
necessary for the implementation of multi-adaptive time-stepping.

Note that we do not solve the dual problem to compute
stability factors (or stability weights) which is necessary to obtain
a reliable error estimate. Thus, the tolerance controls only the size
of the error modulo the stability factor, which is unknown.

In addition, we also compare the two methods for varying size~$L$ of
the domain~$\Omega$, keeping the same initial conditions but scaling
the number of mesh points according to the length of the domain,
$N = 1000L/5$. As the size of the domain increases, we expect the
relative efficiency of the multi-adaptive method to increase, since
the number of inactive components increases relative to the number of
components located within the reaction front.

In Figure~\ref{fig:reaction,cputime}, we plot the CPU time as function
of the tolerance and number of components (size of domain) for the
$\mathrm{mcG}(1)$ and $\mathrm{cG}(1)$ methods. We also summarize the
results in Table~\ref{tab:reaction,tol} and
Table~\ref{tab:reaction,domain}. As expected, the speedup expressed as
the multi-adaptive efficiency index~$\mu$, that is, the ideal speedup
if the cost per degree of freedom were the same for the multi- and
mono-adaptive methods, is large in all test cases, around a
factor~$100$. The speedup in terms of the total number of time slabs
is also large. Note that in Table~\ref{tab:reaction,tol}, the total
number of time slabs~$M$ remains practically constant as the tolerance
and the error are decreased. The decreased tolerance instead results
in finer local resolution of the reaction front, which is evident from
the increasing multi-adaptive efficiency index.  At the same time, the
mono-adaptive method needs to decrease the time step for all
components and so the relative efficiency of the multi-adaptive method
increases as the tolerance decreases. See also
Figure~\ref{fig:reaction,tols} for a comparison of the multi-adaptive
time steps at two different tolerances.

The situation is slightly different in Table~\ref{tab:reaction,domain},
where the tolerance is kept constant but the size of the domain and
number of components vary. Here, the number of time slabs remains
practically constant for both methods, but the multi-adaptive
efficiency index increases as the size of the domain increases, since
the reaction front then becomes more and more localized relative to
the size of the domain. As a result, the efficiency index of the
multi-adaptive method increases as the size of the domain is
increased.

In all test cases, the multi-adaptive method is more efficient than
the standard mono-adaptive method also when the CPU time (wall-clock
time) is chosen as a metric for the comparison. In the first set of
test cases with varying tolerance, the actual speedup is about a
factor $2.0$ whereas in the second test case with varying size of the
domain, the speedup increases from about a factor $2.0$ to a factor
$5.7$ for the range of test cases. These are significant speedups,
although far from the ideal speedup which is given by the
multi-adaptive efficiency index.

There are mainly two reasons that make it difficult to attain full
speedup. The first reason is that as the size of the time slab
increases, the number of iterations~$n$ needed to solve the system of
discrete equations increases. In Table~\ref{tab:reaction,domain}, the
number of iterations, including local iterations on individual
elements as part of a global iteration on the time slab, is about a
factor~$1.5$ larger for the multi-adaptive method. However, the main
overhead lies in the more straightforward implementation of the
mono-adaptive method compared to the more complicated data structures
needed to store and interpolate the multi-adaptive solution. For
constant time step and equal time step for all components, this
overhead is roughly a factor~$5$ for the test problem, but the
overhead increases to about a factor~$100$ when the time slab is
locally refined. It thus remains important to further reduce the
overhead of the implementation in order to increase the range of
problems where the multi-adaptive methods give a positive speedup.

\begin{figure}[htbp]
  \begin{center}
    \psfrag{cputime}{CPU time}
    \psfrag{e}{\hspace{-0.3cm}$\|e(T)\|_{\infty}$}
    \psfrag{N}{$N$}
    \includegraphics[width=6.2cm]{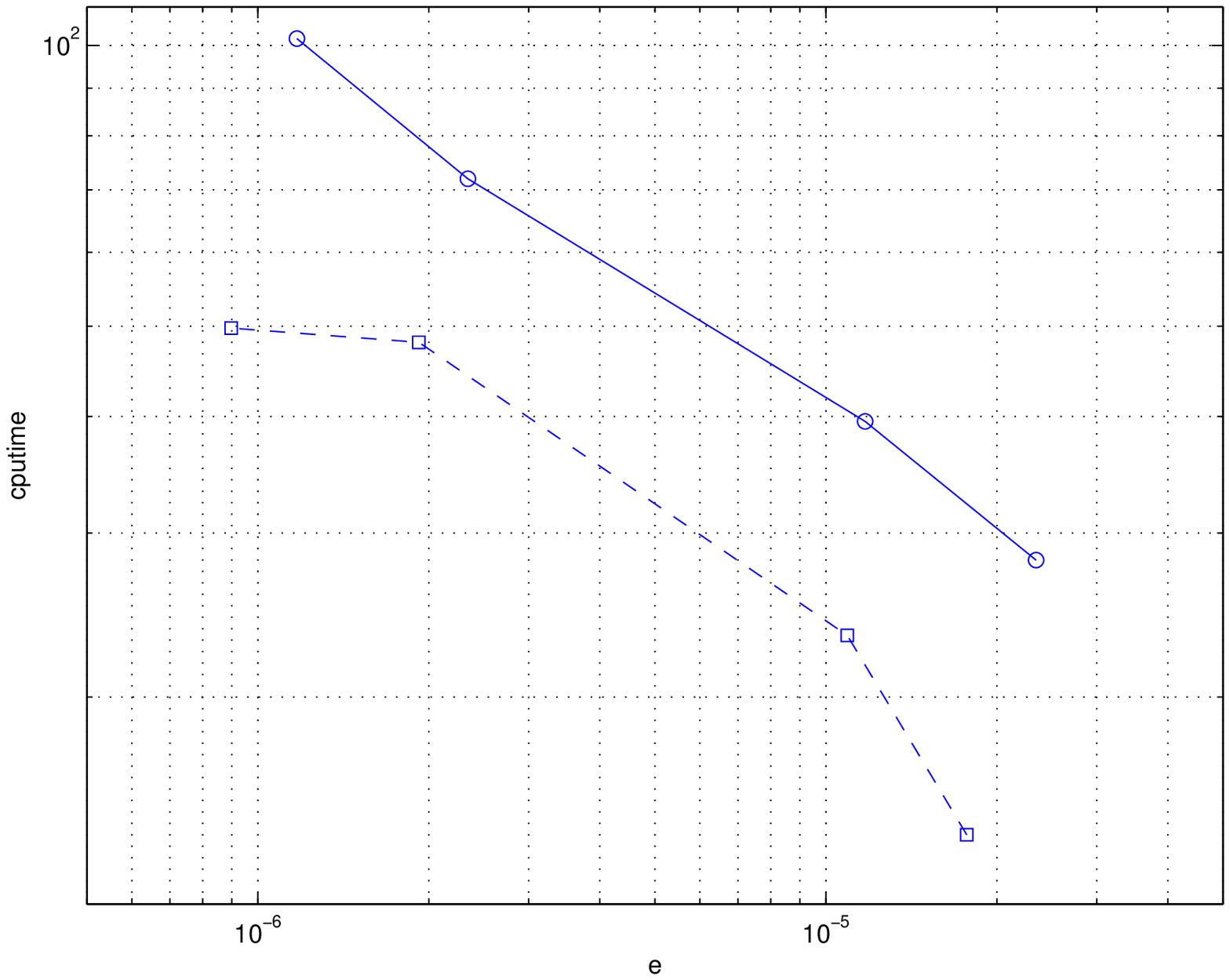}
    \includegraphics[width=6.2cm]{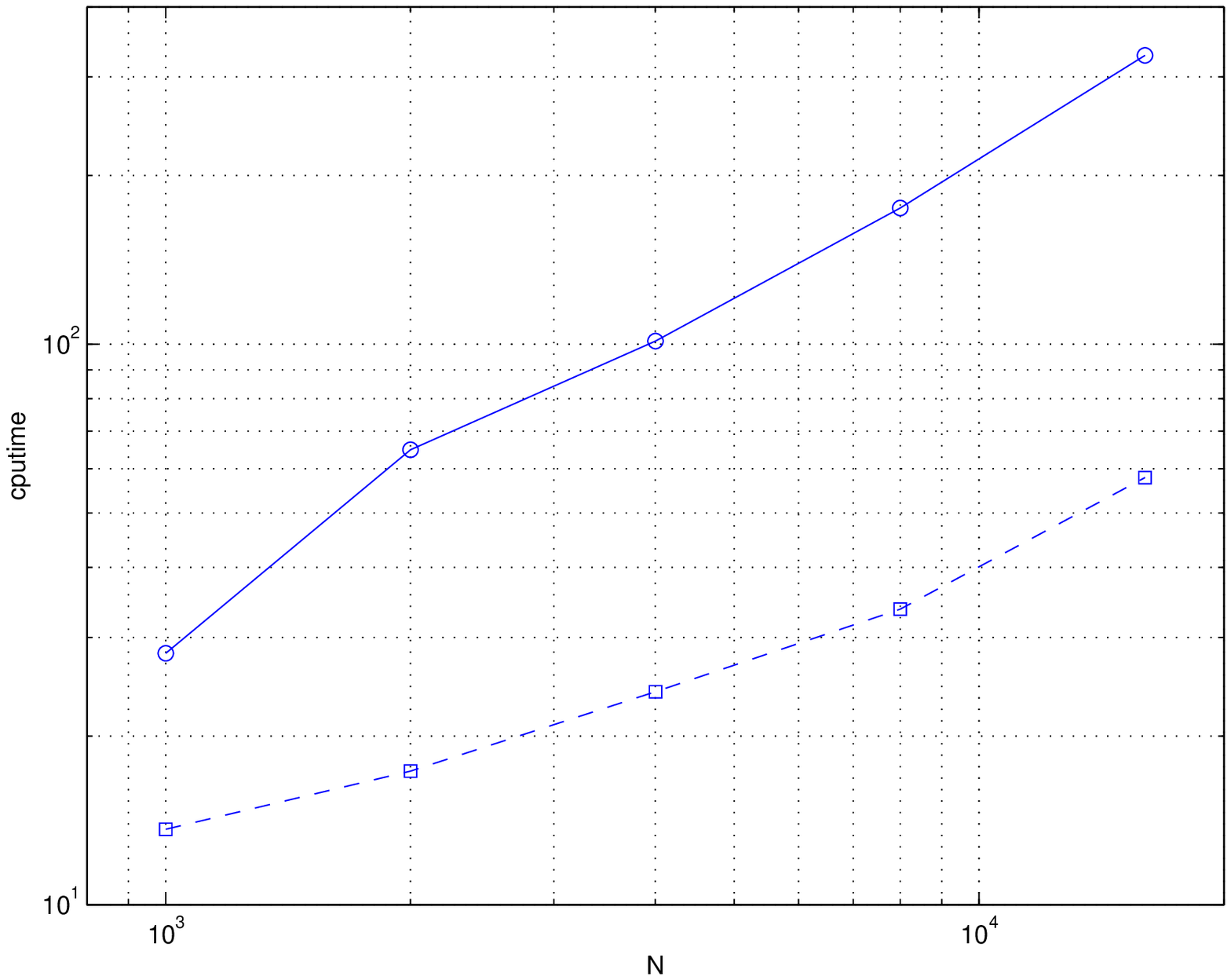}
    \caption{CPU time as function of the error
      (left) and number of components~$N$ (right) for
      $\mathrm{mcG}(1)$ (dashed line) and $\mathrm{cG}(1)$ (solid
      line) for the test problem~(\ref{eq:reaction}).}
    \label{fig:reaction,cputime}
  \end{center}
\end{figure}

\begin{table}[htbp]
  \begin{center}
    \linespread{1.2}\small
    \begin{tabular}{|c|c|c|c|c|c|}
      \hline
      $\mathrm{TOL}$ & $\|e(T)\|_{\infty}$ & CPU time & $M$ & $n$ & $\mu$ \\
      \hline
      \hline
      $1.0\cdot 10^{-6}$ & $1.8\cdot 10^{-5}$ & $14.2 \, \mathrm{s}$ & $1922 \, (5)$ & $3.990 \, (1.498)$ & $ 95.3$ \\
      $5.0\cdot 10^{-7}$ & $1.1\cdot 10^{-5}$ & $23.3 \, \mathrm{s}$ & $1912 \, (9)$ & $4.822 \, (1.544)$ & $138.2$ \\
      $1.0\cdot 10^{-7}$ & $1.9\cdot 10^{-6}$ & $48.1 \, \mathrm{s}$ & $1929 \, (7)$ & $4.905 \, (1.594)$ & $142.6$ \\
      $5.0\cdot 10^{-8}$ & $9.0\cdot 10^{-7}$ & $49.8 \, \mathrm{s}$ & $1917 \, (7)$ & $4.131 \, (1.680)$ & $172.4$ \\
      \hline
      \hline
      $\mathrm{TOL}$ & $\|e(T)\|_{\infty}$ & time & $M$ & $n$ & $\mu$ \\
      \hline
      $1\cdot 10^{-6}$ & $2.3\cdot 10^{-5}$ & $28.1  \, \mathrm{s}$ & $117089 \, (1)$ & $4.0$ & $1.0$ \\
      $5\cdot 10^{-7}$ & $1.2\cdot 10^{-5}$ & $39.5  \, \mathrm{s}$ & $165586 \, (1)$ & $4.0$ & $1.0$ \\
      $1\cdot 10^{-7}$ & $2.3\cdot 10^{-6}$ & $71.9  \, \mathrm{s}$ & $370254 \, (1)$ & $3.0$ & $1.0$ \\
      $5\cdot 10^{-8}$ & $1.2\cdot 10^{-6}$ & $101.7 \, \mathrm{s}$ & $523615 \, (1)$ & $3.0$ & $1.0$ \\
      \hline
    \end{tabular}
    \linespread{1.0}\normalsize
    \caption{Benchmark results for $\mathrm{mcG}(1)$ (above) and
      $\mathrm{cG}(1)$ (below) for varying tolerance and fixed number of
      components $N = 1000$ for the test
      problem~(\ref{eq:reaction}). The table shows the tolerance
      $\mathrm{TOL}$ used for the computation, the error
      $\|e(T)\|_{\infty}$ in the maximum norm at the final time, the
      time used to compute the solution, the number of time slabs $M$
      (with the number of rejected time slabs in parenthesis), the
      average number of iterations $n$ on the time slab system (with
      the number of local iterations on sub-slabs in parenthesis), and
      the multi-adaptive efficiency index~$\mu$.}
    \label{tab:reaction,tol}
  \end{center}
\end{table}

\begin{figure}[htbp]
  \begin{center}
    \psfrag{tol1}{}
    \psfrag{tol2}{}
    \psfrag{x}{$x$}
    \psfrag{k}{$k(x)$}
    \includegraphics[width=12cm]{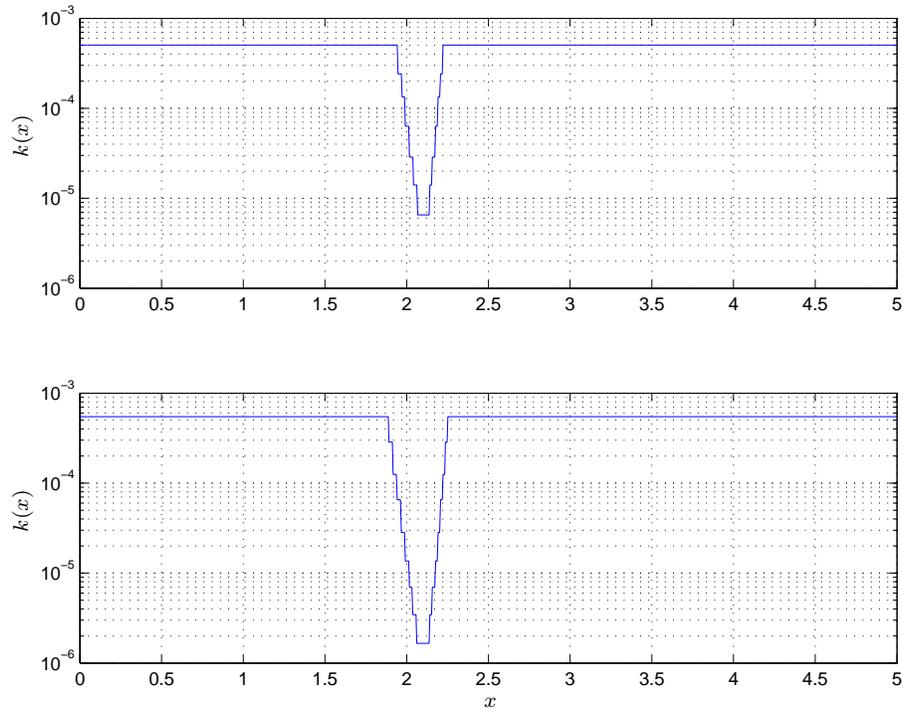}
    \caption{Multi-adaptive time steps at $t = 0.5$ for two different
    tolerances for the test problem~(\ref{eq:reaction}).}
    \label{fig:reaction,tols}
  \end{center}
\end{figure}

\begin{table}[htbp]
  \begin{center}
    \linespread{1.2}\small
    \begin{tabular}{|c|c|c|c|c|c|}
      \hline
      $N$ & $\|e(T)\|_{\infty}$ & CPU time & $M$ & $n$ & $\mu$ \\
      \hline
      \hline
      $1000$  & $1.8\cdot 10^{-5}$ & $13.6 \, \mathrm{ s}$ & $1922 \, (5)$ & $4.0 \, (1.5)$ & $95.3$  \\
      $2000$  & $1.7\cdot 10^{-5}$ & $17.3 \, \mathrm{ s}$ & $1923 \, (5)$ & $4.0 \, (1.2)$ & $140.5$ \\
      $4000$  & $1.6\cdot 10^{-5}$ & $24.0 \, \mathrm{ s}$ & $1920 \, (6)$ & $4.0 \, (1.0)$ & $185.0$ \\
      $8000$  & $1.7\cdot 10^{-5}$ & $33.7 \, \mathrm{ s}$ & $1918 \, (5)$ & $4.0 \, (1.0)$ & $218.8$ \\
      $16000$ & $1.7\cdot 10^{-5}$ & $57.9 \, \mathrm{ s}$ & $1919 \, (5)$ & $4.0 \, (1.0)$ & $240.0$ \\
      \hline
      \hline
      $N$ & $\|e(T)\|_{\infty}$ & time & $M$ & $n$ & $\mu$ \\
      \hline
      $1000$  & $2.3\cdot 10^{-5}$ & $28.1  \, \mathrm{s}$ & $117089 \, (1)$ & $4.0$ & $ 1.0$ \\
      $2000$  & $2.2\cdot 10^{-5}$ & $64.8  \, \mathrm{s}$ & $117091 \, (1)$ & $4.0$ & $ 1.0$ \\
      $4000$  & $2.2\cdot 10^{-5}$ & $101.3 \, \mathrm{s}$ & $117090 \, (1)$ & $4.0$ & $ 1.0$ \\
      $8000$  & $2.2\cdot 10^{-5}$ & $175.1 \, \mathrm{s}$ & $117089 \, (1)$ & $4.0$ & $ 1.0$ \\
      $16000$ & $2.2\cdot 10^{-5}$ & $327.7 \, \mathrm{s}$ & $117089 \, (1)$ & $4.0$ & $ 1.0$ \\
      \hline
    \end{tabular}
    \linespread{1.0}\normalsize
    \caption{Benchmark results for $\mathrm{mcG}(1)$ (above) and $\mathrm{cG}(1)$ (below) for fixed tolerance
      $\mathrm{TOL} = 1.0\cdot{}10^{-6}$ and varying number of
      components (and size of domain). (See
      Table~\ref{tab:reaction,tol} for an explanation of table legends.)}
    \label{tab:reaction,domain}
  \end{center}
\end{table}

\subsection{The wave equation}

Next, we consider the wave equation,
\begin{equation} \label{eq:wave}
  \begin{split}
    u_{tt} - \Delta u &= 0 \quad \mbox{ in } \Omega \times (0,T], \\
    \partial_n u &= 0 \quad \mbox{ on } \partial \Omega \times (0,T], \\
    u(\cdot, 0) &= u_0 \quad \mbox{ in } \Omega,
  \end{split}
\end{equation}
on a two-dimensional domain~$\Omega$ consisting of two square sub-domains of
side length~$0.5$ separated by a thin wall with a narrow slit of size
$0.0001 \times 0.0001$ at its center. The initial condition is chosen
as a plane wave traversing the domain from right to left.
In Figure~\ref{fig:wave,u,k}, we plot the initial data together with
the (fixed) multi-adaptive time steps.
The resulting solution is shown in Figure~\ref{fig:wave,solution}.

The geometry of the domain $\Omega$ forces the discretization to be
very fine close to the narrow slit. Further away from the slit, we let
the mesh be coarse. The mesh was created by specifying a mesh size $h$
with $h \gg w$ where $w$ is the width of the narrow
slit. We note that for the multi-adaptive efficiency index $\mu$
defined in~(\ref{eq:efficiencyindex}) to be large, the total number of
elements must be large in comparison the to number of small elements
close to the narrow slit. Furthermore, the average mesh size must be
large compared to the mesh size close to the narrow slit.

For a mono-adaptive method, a global CFL condition puts a limit on the
size of the global time step, roughly given by
\begin{equation}
  k \leq h_{\min} = \min_{x\in\Omega} h(x),
\end{equation}
where $h = h(x)$ is the local mesh size. With a larger time step, an
explicit method will be unstable or, correspondingly, direct
fixed-point iteration on the system of discrete equations on each
time slab will not converge without suitable stabilization.

On the other hand, with a multi-adaptive method, the time step may be
chosen to satisfy the CFL condition only locally, that is,
\begin{equation}
  k(x) \leq h(x), \quad x \in \Omega,
\end{equation}
and as a result, the number of local steps may decrease significantly
(depending on the properties of the mesh). In this case, with $k = 0.1
h$, the speedup for the multi-adaptive~$\mathrm{mcG}(1)$ method was a
factor~$4.2$.

\begin{figure}[htbp]
  \begin{center}
    \includegraphics[width=6.2cm]{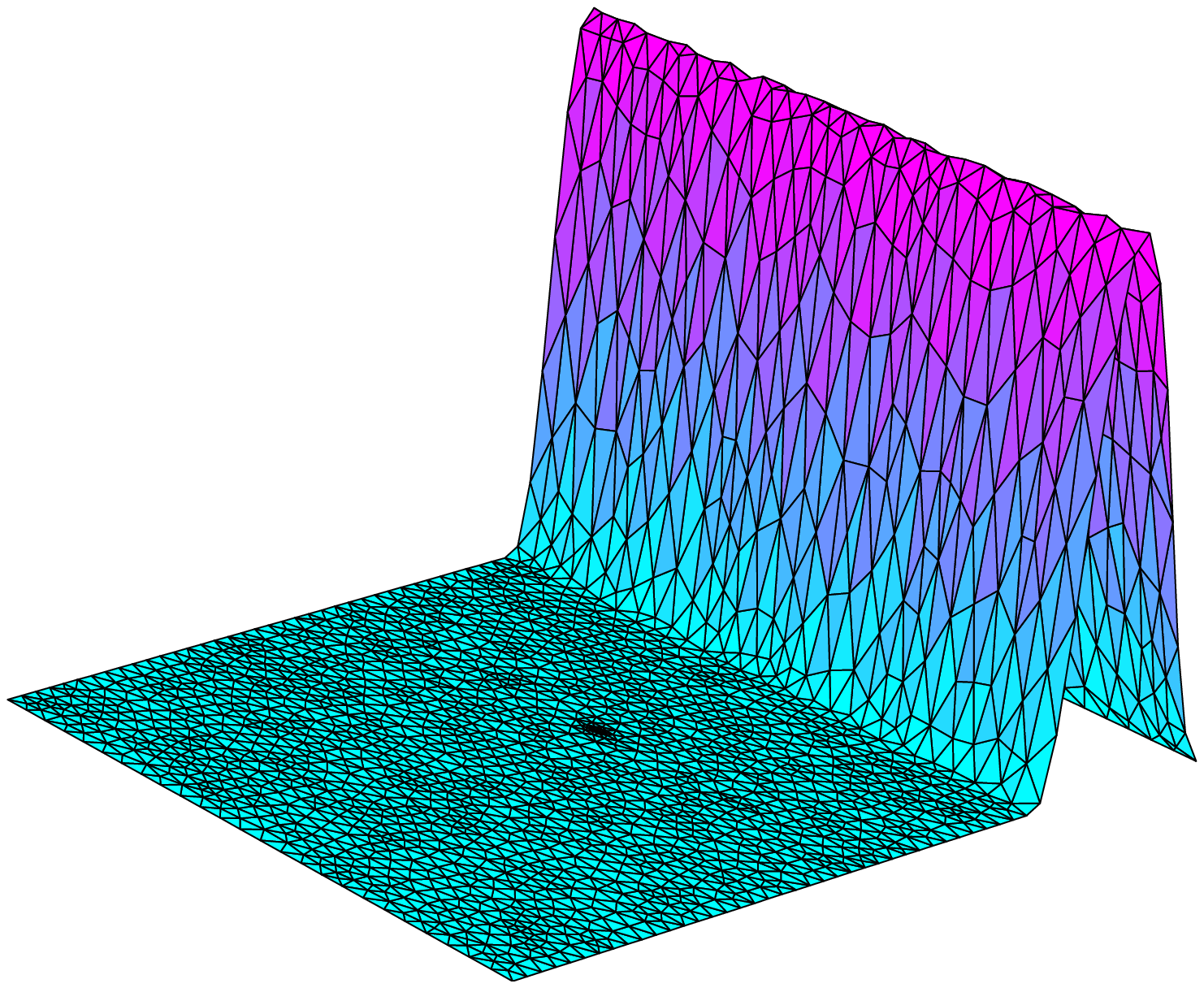}
    \includegraphics[width=6.2cm]{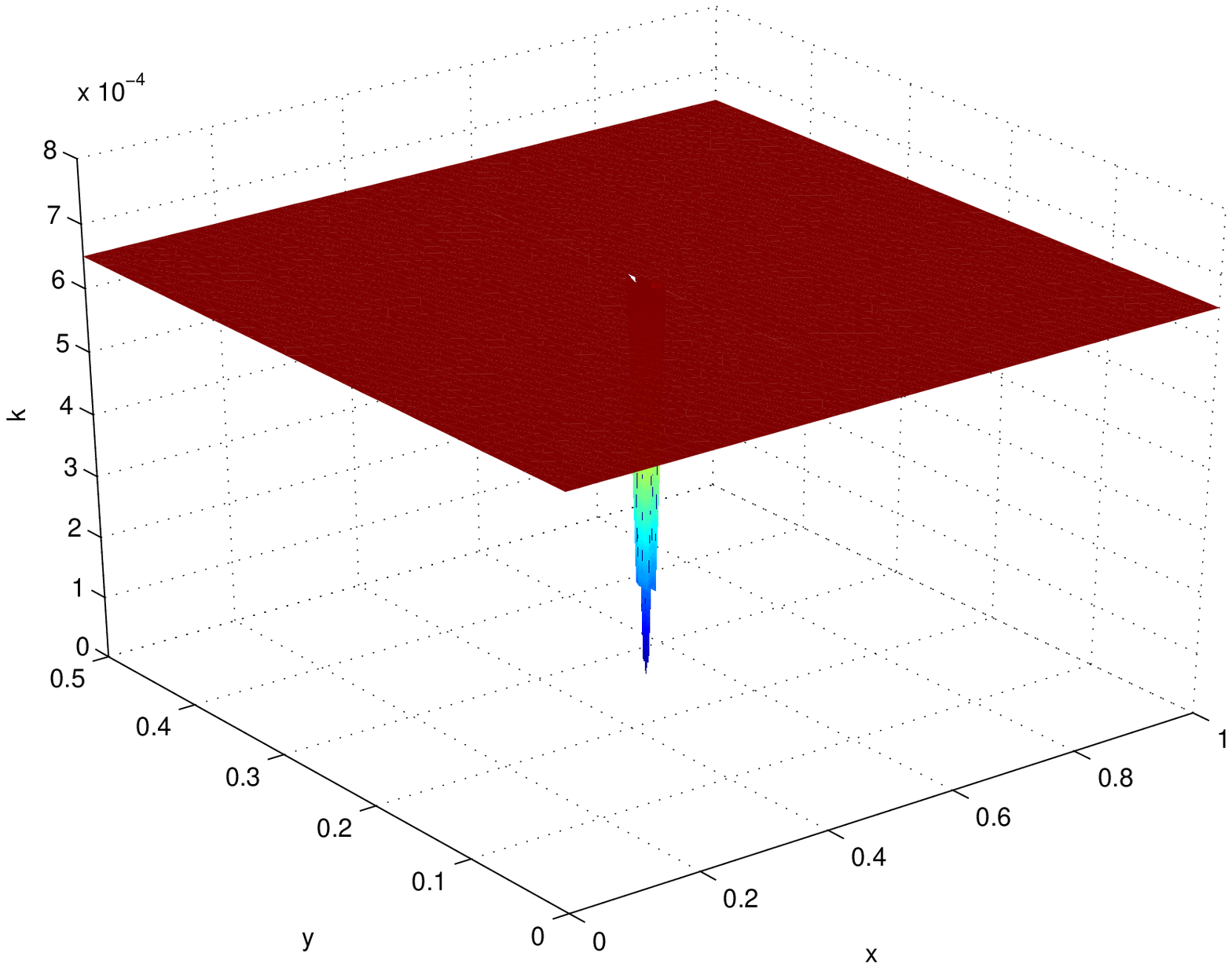}
    \caption{Initial data (left) and multi-adaptive time steps (right)
      for the solution of the wave equation.}
    \label{fig:wave,u,k}
  \end{center}
\end{figure}

\begin{figure}[htbp]
  \begin{center}
    \includegraphics[width=6.2cm]{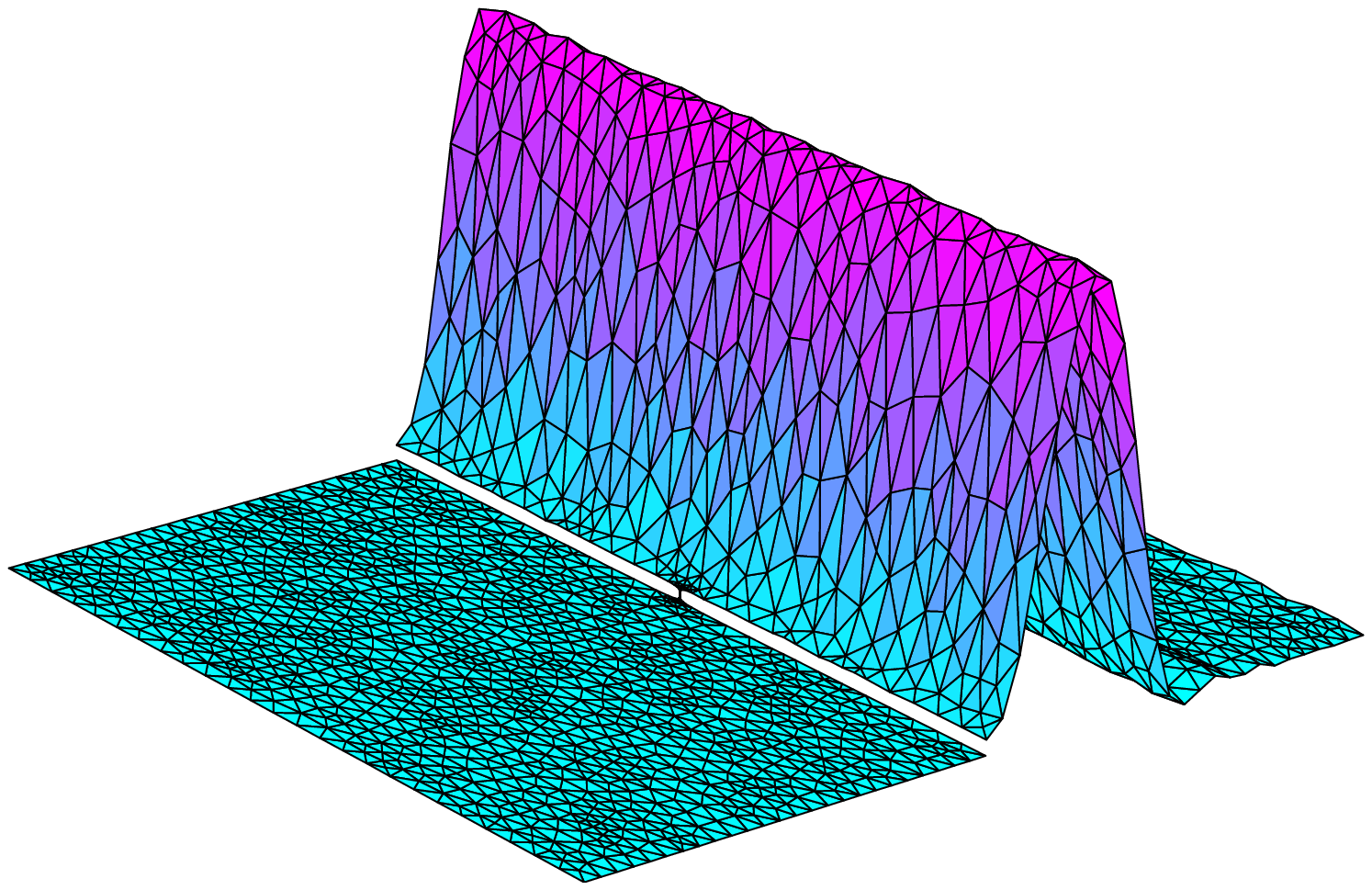}
    \includegraphics[width=6.2cm]{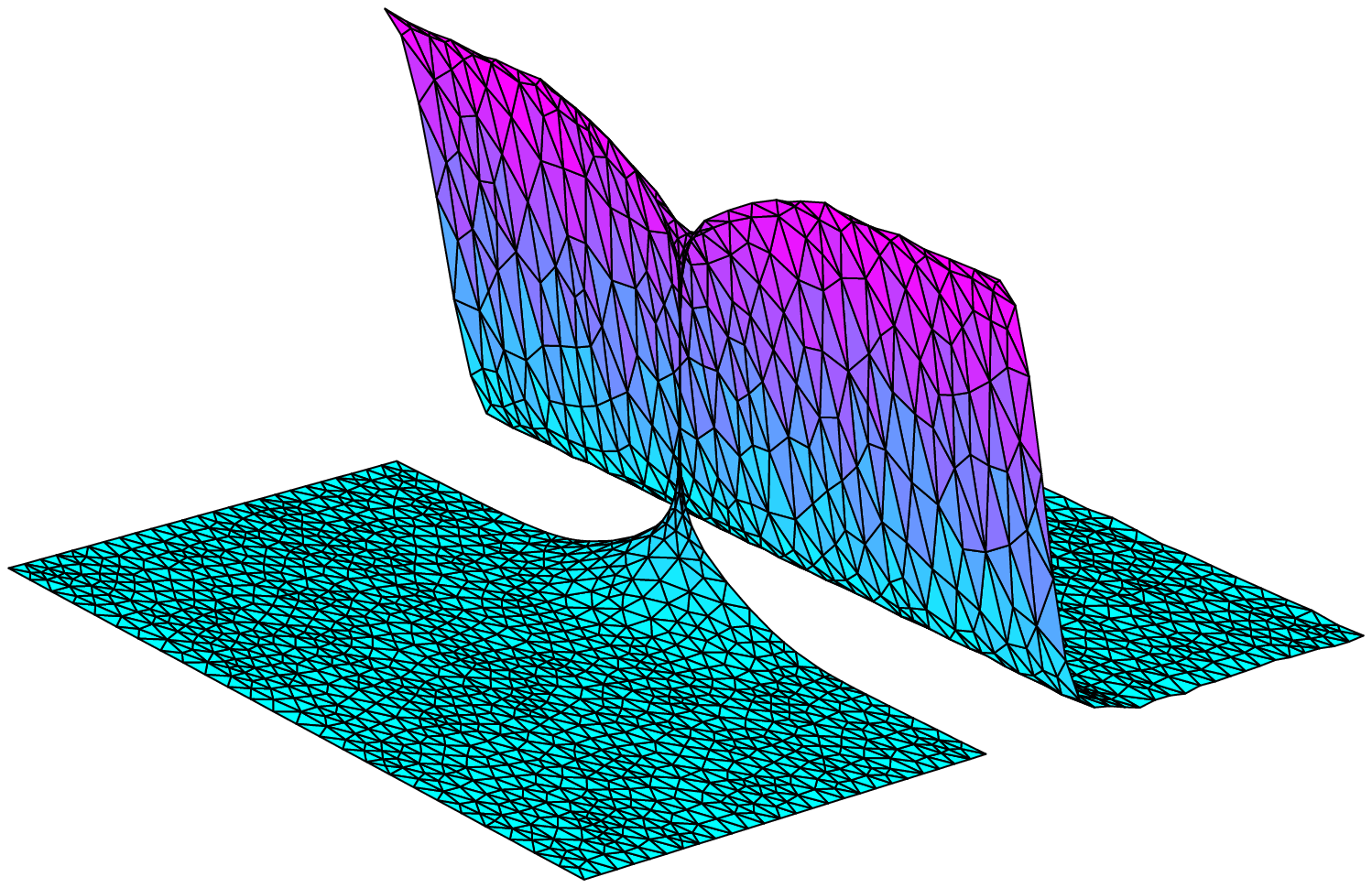}\\
    \includegraphics[width=6.2cm]{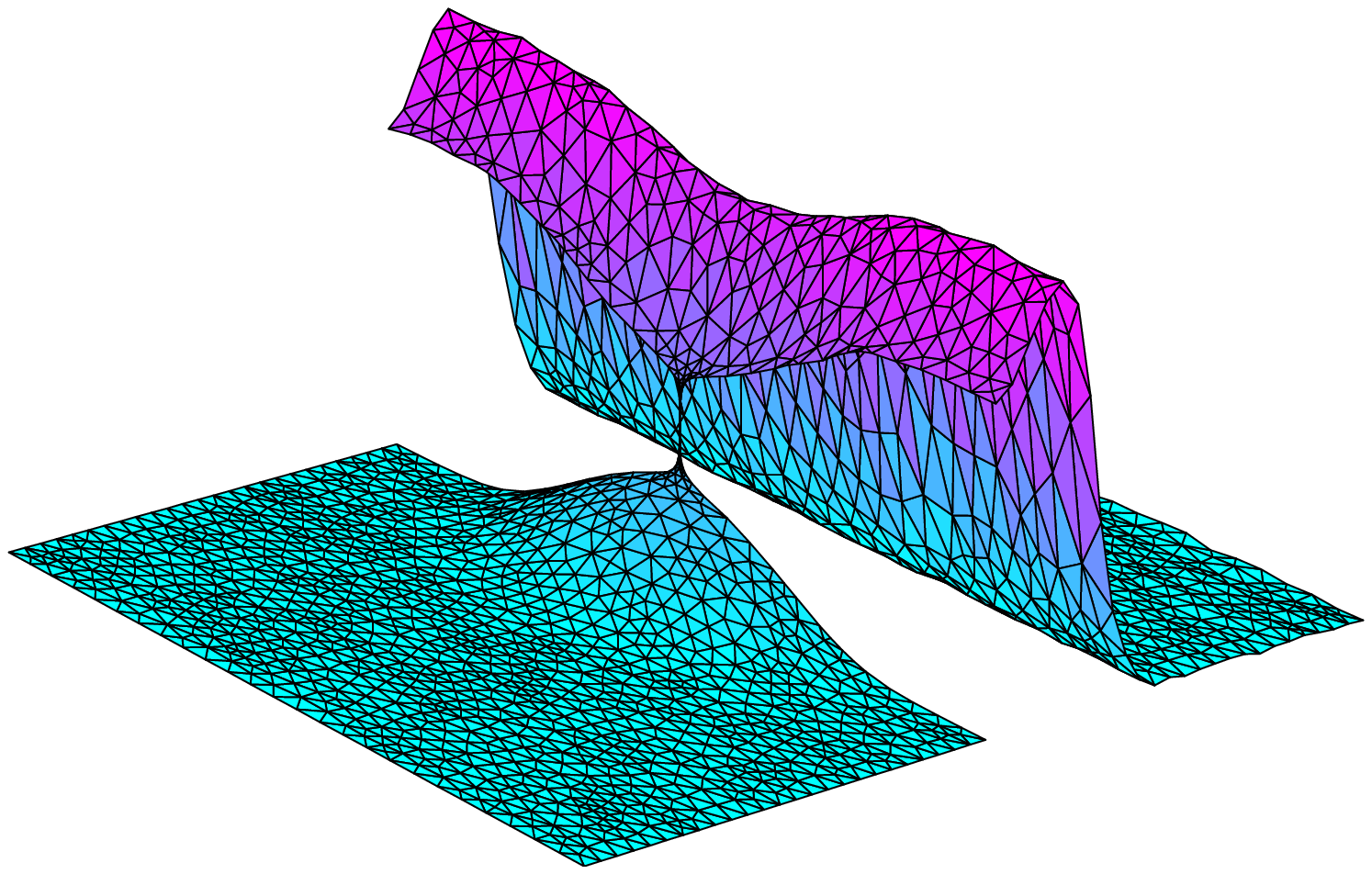}
    \includegraphics[width=6.2cm]{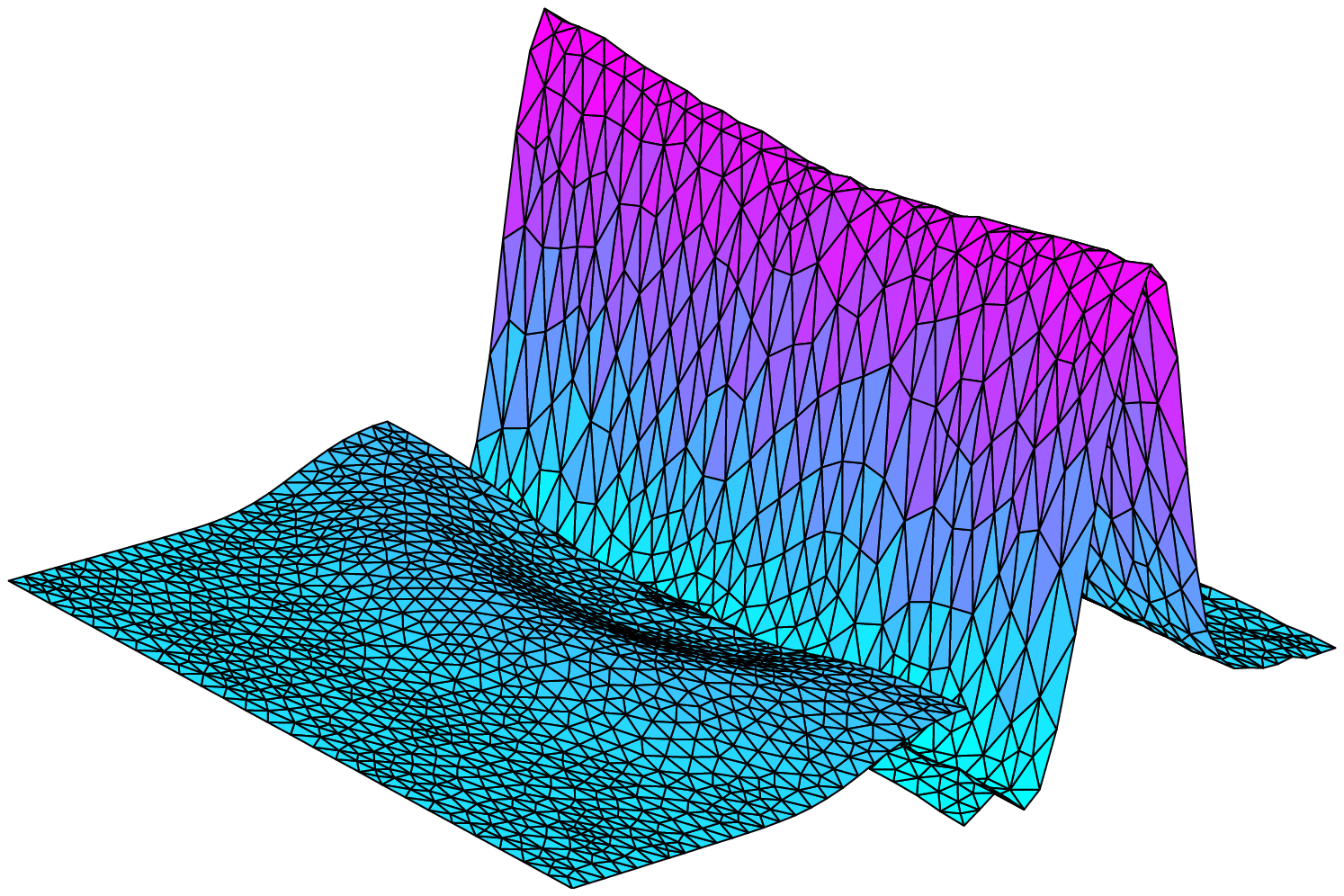}
    \caption{The solution of the wave equation at times $t = 0.25$,
      $t = 0.4$, $t = 0.45$ and $t = 0.6$.}
    \label{fig:wave,solution}
  \end{center}
\end{figure}

\begin{figure}[htbp]
  \begin{center}
    \psfrag{x}{}
    \psfrag{y}{}
    \includegraphics[width=6.2cm]{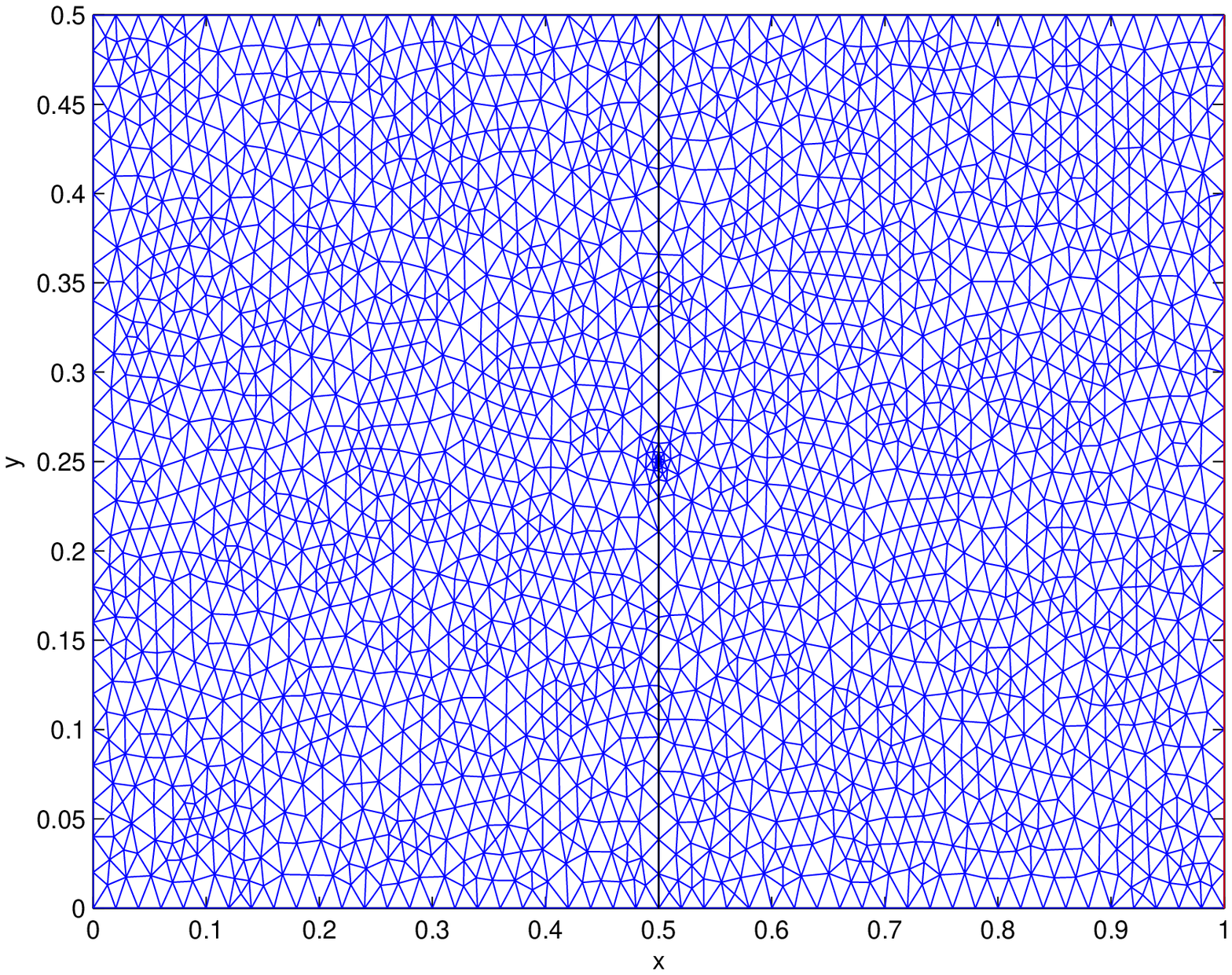}
    \includegraphics[width=6.2cm]{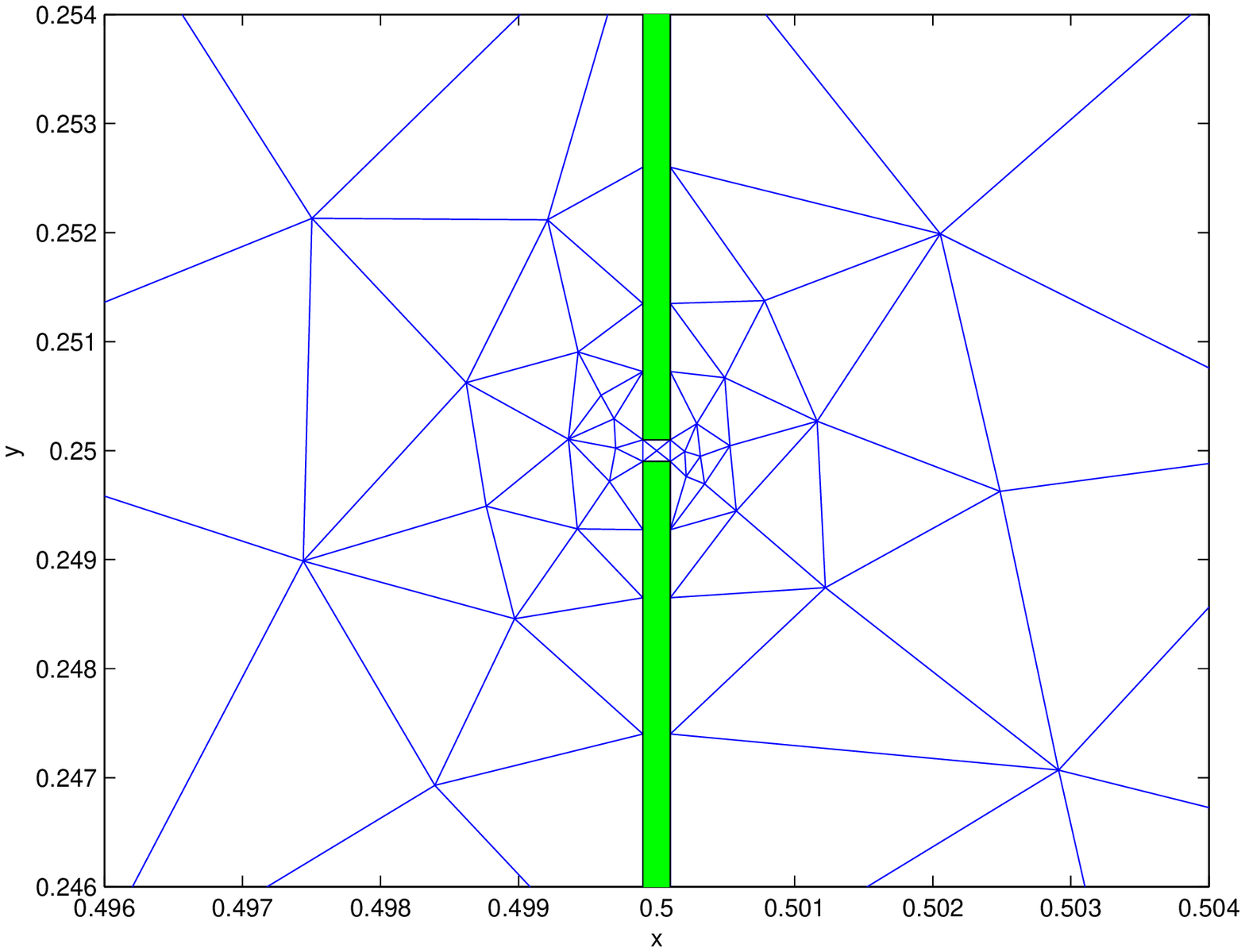}
    \caption{The mesh used for the solution of the wave equation on a
      domain intersected by a thin wall with a narrow slit (left) and
      details of the mesh close to the slit (right).}
    \label{fig:wave,mesh}
  \end{center}
\end{figure}

\section{Conclusions}
\label{sec:conclusions}

We have presented algorithms and data structures for multi-adaptive
time-stepping, including the recursive construction of time slabs and
efficient interpolation of multi-adaptive solutions. The efficiency of
the multi-adaptive methods was demonstrated for a pair of benchmark
problems. The multi-adaptive methods \mcgq{} and \mdgq{} are available
as components of DOLFIN, together with implementations of the standard
mono-adaptive \cgq{} and \dgq{} methods. The ODE solvers of DOLFIN are
currently being integrated with other components of the FEniCS
project, in particular the FEniCS Form Compiler
(FFC)~\cite{logg:www:04,KirLog2006,KirLog2007} in order to provide
reliable, efficient and automatic integration of time dependent PDEs.

\bibliographystyle{acmtrans}
\bibliography{bibliography}

\end{document}